\input epsf
\magnification=\magstep1
\documentstyle{amsppt}
\pagewidth{6.48truein}
\pageheight{9truein}
\TagsOnRight
\NoRunningHeads
\catcode`\@=11
\def\logo@{}
\footline={\ifnum\pageno>1 \hfil\folio\hfil\else\hfil\fi}
\topmatter
\title A generalization of Kuo condensation
\endtitle
\vskip-0.2in
\author Mihai Ciucu\endauthor
\thanks Research supported in part by NSF grant DMS-1101670.
\endthanks
\affil
  Department of Mathematics, Indiana University\\
  Bloomington, Indiana 47405
\endaffil
\abstract
Kuo introduced his 4-point condensation in 2003 for bipartite planar graphs. In 2006 Kuo generalized this 4-point condensation to planar graphs that are not necessarily bipartite. His formula expressed the product between the number of perfect matching of the original graph $G$ and that of the subgraph obtained from $G$ by removing the four distinguished vertices as a Pfaffian of order 4, whose entries are numbers of perfect matchings of subgraphs of $G$ obtained by removing various pairs of vertices chosen from among the four distinguished ones. The compelling elegance of this formula is inviting of generalization. Kuo generalized it to $2k$ points under the special assumption that the subgraph obtained by removing some subset of the $2k$ vertices has precisely one perfect matching. In this paper we prove that the formula holds in the general case. We also present a number of applications.
\endabstract
\endtopmatter

\document

\def\mysec#1{\bigskip\centerline{\bf #1}\message{ * }\nopagebreak\bigskip\par}

\def\myref#1{\item"{[{\bf #1}]}"}

\def\pf{{\it Proof.\ }}

\def\epf{\hfill{$\square$}\smallpagebreak}

\def\cite#1{\relaxnext@
  \def\nextiii@##1,##2\end@{[{\bf##1},\,##2]}%
  \in@,{#1}\ifin@\def\next{\nextiii@#1\end@}\else
  \def\next{[{\bf#1}]}\fi\next}
\def\proclaimheadfont@{\smc}

\def\pf{{\it Proof.\ }}

\define\C{{\Bbb C}}
\define\Rep{\operatorname{Re}}

\define\M{\operatorname{M}}

\define\Pf{\operatorname{Pf}}


\define\twoline#1#2{\line{\hfill{\smc #1}\hfill{\smc #2}\hfill}}
\define\twolinetwo#1#2{\line{{\smc #1}\hfill{\smc #2}}}
\define\twolinethree#1#2{\line{\phantom{poco}{\smc #1}\hfill{\smc #2}\phantom{poco}}}
\define\threeline#1#2#3{\line{\hfill{\smc #1}\hfill{\smc #2}\hfill{\smc #3}\hfill}}

\def\mypic#1{\epsffile{#1}}



\define\anglep{1}
\define\ff{2}
\define\hexnotch{3}
\define\CLP{4}
\define\Eisen{5}
\define\Fulmek{6}
\define\GT{7}
\define\KuoOne{8}
\define\KuoTwo{9}
\define\MacM{10}
\define\Sl{11}

\define\eaa{1.1}
\define\eab{1.2}
\define\eac{1.3}

\define\eba{2.1}
\define\ebb{2.2}
\define\ebc{2.3}
\define\ebd{2.4}
\define\ebe{2.5}
\define\ebf{2.6}
\define\ebg{2.7}
\define\ebh{2.8}
\define\ebi{2.9}
\define\ebj{2.10}

\define\eca{3.1}
\define\ecb{3.2}
\define\ecc{3.3}
\define\ecdd{3.4}
\define\ecddd{3.5}
\define\ecd{3.6}
\define\ece{3.7}
\define\ecf{3.8}

\define\ech{3.10}

\define\eda{4.1}
\define\edc{4.2}
\define\edd{4.3}
\define\ede{4.4}
\define\edf{4.5}
\define\edg{4.6}
\define\edh{4.7}
\define\edi{4.8}
\define\edj{4.9}


\define\tba{2.1}
\define\tbb{2.2}
\define\tbc{2.3}
\define\tbd{2.4}

\define\tca{3.1}
\define\tcb{3.2}

\define\tda{4.1}
\define\tdb{4.2}
\define\tdc{4.3}
\define\tdd{4.4}


\define\fba{2.1}

\define\fca{3.1}
\define\fcb{3.2}
\define\fcc{3.3}
\define\fcd{3.4}
\define\fce{3.5}
\define\fcf{3.6}
\define\fcg{3.7}

\define\fda{4.1}
\define\fdb{4.2}
\define\fdc{4.3}
\define\fdd{4.4}
\define\fddp{4.5}
\define\fde{4.6}
\define\fddpp{4.7}
\define\fddppp{4.8}

\vskip-0.05in
\mysec{1. Introduction}

In \cite{\KuoOne} Kuo introduced the method of graphical condensation as a powerful way to obtain recurrences for the number of perfect matchings of planar bipartite graphs. Let $G$ be a plane bipartite graph with the same number of vertices in its two color classes $V_1$ and $V_2$. Let $a$, $b$, $c$, $d$ be vertices appearing in cyclic order on some face of $G$, with $a,c\in V_1$ and $b,d\in V_2$. Then \cite{\KuoOne,Theorem\,2.1} states that
$$
\M(G)\M(G\setminus\{a,b,c,d\})=
\M(G\setminus\{a,b\})\M(G\setminus\{c,d\})+
\M(G\setminus\{a,d\})\M(G\setminus\{b,c\}),
\tag\eaa
$$
where $\M(H)$ stands for the number of perfect matchings of the graph $H$.

Kuo then generalized this in \cite{\KuoTwo} to planar graphs that are not necessarily bipartite. Namely, for any planar graph $G$ and any four vertices $a$, $b$, $c$, $d$ that appear in cyclic order on some face of $G$, one has by \cite{\KuoTwo,Proposition\,1.1} that
$$
\spreadlines{3\jot}
\align
\M(G)\M(G\setminus\{a,b,c,d\})
=
\M(G\setminus\{a,b\})\M(G\setminus\{c,d\})
&
-
\M(G\setminus\{a,c\})\M(G\setminus\{b,d\})
\\
&
+
\M(G\setminus\{a,d\})\M(G\setminus\{b,c\}),
\tag\eab
\endalign
$$
which, as Kuo points out in \cite{\KuoTwo}, can also be written in the compelling form
$$
\spreadlines{3\jot}
\align
&
\M(G)\M(G\setminus\{a,b,c,d\})
=
\\
&\ \ \ \ \ \ \ \ 
\Pf\left[\matrix
0 & \M(G\setminus\{a,b\}) & \M(G\setminus\{a,c\}) & \M(G\setminus\{a,d\})\\
- \M(G\setminus\{a,b\}) & 0 & \M(G\setminus\{b,c\}) & \M(G\setminus\{b,d\})\\
-\M(G\setminus\{a,c\}) & -\M(G\setminus\{b,c\}) & 0 & \M(G\setminus\{c,d\})\\
-\M(G\setminus\{a,d\}) & -\M(G\setminus\{b,d\}) & -\M(G\setminus\{c,d\}) & 0
\endmatrix\right].
\tag\eac
\endalign
$$

The striking elegance of this formula is inviting of generalization. Kuo generalized it to $2k$ points under the special assumption that the subgraph obtained by removing some subset of the $2k$ vertices has precisely one perfect matching (see \cite{\KuoTwo, Theorem\,3.1}).

In this paper we prove that the formula holds in the general case. As applications of it, we present a conceptual proof of a theorem of Eisenk\"olbl and a generalization of it. For three recent applications of Kuo's original formula see \cite{\anglep}, \cite{\ff} and \cite{\hexnotch}.

\mysec{2. The general Pfaffian graphical condensation}

Our generalization of Kuo's graphical condensation (\eac) is the following. A {\it weighted} graph is a graph with weights (that could be considered indeterminates) on its edges. For a weighted graph $G$, $\M(G)$ denotes the sum of the weights of the perfect matchings of $G$, where the weight of a perfect matching is taken to be the product of the weights of its constituent edges.

\proclaim{Theorem \tba} Let $G$ be a planar graph with the vertices $a_1,\dotsc,a_{2k}$ appearing in that cyclic order on a face of $G$. Consider the skew-symmetric matrix $A=(a_{ij})_{1\leq i,j\leq 2k}$ with entries given by
$$
a_{ij}:=\left\{\matrix
\M(G\setminus\{a_i,a_j\}),\ \ \ \text{\rm if $i<j$},\\
-\M(G\setminus\{a_i,a_j\}),\ \ \ \text{\rm if $i>j$}.
\endmatrix\right.
$$
Then we have that
$$
\left[\M(G)\right]^{k-1}\M(G\setminus\{a_1,\dotsc,a_{2k}\})
=
\Pf(A).
\tag\eba
$$

\endproclaim

In our proof of the above theorem we make use of the following auxiliary result that presents some interest on its own. 

\proclaim{Proposition \tbb} Let $G$ be a planar graph with the vertices $a_1,\dotsc,a_{2k}$ appearing in that cyclic order on a face of $G$. Then
$$
\spreadlines{3\jot}
\align
&\!\!\!\!\!\!\!\!\!\!\!\!\!\!\!\!
\M(G)\M(G\setminus\{a_1,\dotsc,a_{2k}\})
+
\M(G\setminus\{a_1,a_3\})\M(G\setminus\overline{\{a_1,a_3\}})
+
\cdots
\\
&\ \ \ \ \ \ \ \ \ \ \ \ \ \ \ \ \ \ \ \ \ \ \ \ \ \ \ \ \ \ \ \ \ \ \ \ \ \ \ \ \ \ 
+\M(G\setminus\{a_1,a_{2k-1}\})\M(G\setminus\overline{\{a_1,a_{2k-1}\}})
\\
&
=
\M(G\setminus\{a_1,a_2\})\M(G\setminus\overline{\{a_1,a_2\}})
+
\M(G\setminus\{a_1,a_4\})\M(G\setminus\overline{\{a_1,a_4\}})
+
\cdots
\\
&\ \ \ \ \ \ \ \ \ \ \ \ \ \ \ \ \ \ \ \ \ \ \ \ \ \ \ \ \ \ \ \ \ \ \ \ \ \ \ \ \ \ 
+\M(G\setminus\{a_1,a_{2k}\})\M(G\setminus\overline{\{a_1,a_{2k}\}}),
\tag\ebb
\endalign
$$
where $\overline{\{a_i,a_j\}}$ stands for the complement of $\{a_i,a_j\}$ in the set $\{a_1,\dotsc,a_{2k}\}$.

\endproclaim

\pf Denote by $\Cal M(G)$ the set of perfect matchings of the graph $G$. Patterned on the two sides of equation (\ebb), consider the disjoint unions of Cartesian products
$$
\spreadlines{3\jot}
\align
&\!\!\!\!\!\!\!\!\!\!\!\!\!\!\!\!
\Cal M(G)\times \Cal M(G\setminus\{a_1,\dotsc,a_{2k}\}) 
\cup
\Cal M(G\setminus\{a_1,a_3\})\times\Cal M(G\setminus\overline{\{a_1,a_3\}})
\cup
\cdots
\\
&\ \ \ \ \ \ \ \ \ \ \ \ \ \ \ \ \ \ \ \ \ \ \ \ \ \ \ \ \ \ \ \ \ \ \ \ \ \ \ \ \ \ 
\cup
\Cal M(G\setminus\{a_1,a_{2k-1}\})\times\Cal M(G\setminus\overline{\{a_1,a_{2k-1}\}})
\tag\ebc
\endalign
$$
and
$$
\spreadlines{3\jot}
\align
&\!\!\!\!\!\!\!\!\!\!\!\!\!\!\!\!
\Cal M(G\setminus\{a_1,a_2\})\times\Cal M(G\setminus\overline{\{a_1,a_2\}})
\cup
\Cal M(G\setminus\{a_1,a_4\})\times\Cal M(G\setminus\overline{\{a_1,a_4\}})
\cup
\cdots
\\
&\ \ \ \ \ \ \ \ \ \ \ \ \ \ \ \ \ \ \ \ \ \ \ \ \ \ \ \ \ \ \ \ \ \ \ \ \ \ \ \ \ \ 
\cup
\Cal M(G\setminus\{a_1,a_{2k}\})\times\Cal M(G\setminus\overline{\{a_1,a_{2k}\}}).
\tag\ebd
\endalign
$$
For any element $(\mu,\nu)$ of (\ebc) or (\ebd), think of the edges of $\mu$ as being marked by solid lines, and of the edges of $\nu$ as marked by dotted lines, {\it on the same copy of the graph $G$} (any edge common to $\mu$ and $\nu$ will be marked both solid and dotted, by two parallel arcs). 

Define the weight of $(\mu,\nu)$ to be the product of the weight of $\mu$ and the weight of $\nu$. Then the total weight of the elements of the set (\ebc) is equal to the left hand side of equation (\ebb), while the total weight of the elements of the set (\ebd) equals the right hand side of (\ebb). Therefore, to prove (\ebb) it suffices to construct a weight-preserving bijection between the sets (\ebc) and (\ebd).

We construct such a bijection as follows. Let $(\mu,\nu)$ be an element of (\ebc). Our construction depends upon the particular set of the union (\ebc) that $(\mu,\nu)$ belongs to.

If $(\mu,\nu)\in\Cal M(G)\times \Cal M(G\setminus\{a_1,\dotsc,a_{2k}\})$, map  $(\mu,\nu)$ to what we get from it by ``shifting along the path containing $a_1$.'' More precisely, note that when considering the edges of $\mu$ and $\nu$ together on the same copy of $G$, each of the vertices $a_1,\dotsc,a_{2k}$ is incident to precisely one edge (namely, a solid edge), while all the other vertices of $G$ are incident to one solid edge and one dotted edge. This implies that $\mu\cup\nu$ is the disjoint union of paths connecting the $a_i$'s to one another in pairs, and cycles covering the remaining vertices of~$G$. Consider the path containing $a_1$, and change each solid edge in it to dotted, and each dotted edge to solid. Denote the resulting pair of matchings by $(\mu',\nu')$. 

Since before the reversal of colors the end edges of this path were solid, after the reversal they are both dotted. In addition, this path must connect $a_1$ to one of $a_2,a_4,\dotsc,a_{2k}$, because if it connected $a_1$ to an odd-indexed $a_{2i+1}$ that would isolate the $2i-1$ vertices $a_2,a_3,\dotsc,a_{2i}$ from the other $a_j$'s, making it impossible for them to be connected up by disjoint paths. Therefore, $(\mu',\nu')$ is an element of (\ebd).

Suppose now that $(\mu,\nu)\in\Cal M(G\setminus\{a_1,a_3\})\times \Cal M(G\setminus\overline{\{a_1,a_3\}})$. Then we map $(\mu,\nu)$ to the pair of matchings $(\mu',\nu')$ obtained from it by reversing ``solid'' and ``dotted'' along the path $P$ in $\mu\cup\nu$ containing $a_3$. By the argument in the previous paragraph, this path must connect $a_3$ to one of $a_2,a_4,\dotsc,a_{2k}$. Note that, before the reversal, the end edge of this path incident to $a_3$ was dotted, and the other end edge was solid. Therefore, after the reversal, the other end point of the path $P$ (which is one of  $a_2,a_4,\dotsc,a_{2k}$) swaps places with $a_3$ from the point of view of being matched by a solid line versus a dotted line, and thus $(\mu',\nu')$ is an element of (\ebd).

Finally, if $(\mu,\nu)\in\Cal M(G\setminus\{a_1,a_{2i+1}\})\times \Cal M(G\setminus\overline{\{a_1,a_{2i+1}\}})$ with $i>1$, use the construction in the previous paragraph with $a_3$ replaced by $a_{2i+1}$.

The map $(\mu,\nu)\mapsto(\mu',\nu')$ described above can easily be inverted. Indeed, given an element $(\mu',\nu')$ of the union (\ebd), the pair $(\mu,\nu)$ that gets mapped to it is obtained by shifting along the path in $\mu'\cup\nu'$ that contains the vertex $a_{2i}$, where $i$ is the index for which $(\mu',\nu')\in\Cal M(G\setminus\{a_1,a_{2i}\})\times\Cal M(G\setminus\overline{\{a_1,a_{2i}\}})$. Since the map $(\mu,\nu)\mapsto(\mu',\nu')$ is also clearly weight-preserving, this completes the proof. \epf

We will also need the following classical Pfaffian analog of the expansion of a determinant along a row.

\proclaim{Lemma \tbc} For any $2n\times2n$ skew-symmetric matrix $A=(a_{ij})$, we have
$$
\Pf(A)=\sum_{i=2}^{2n} (-1)^i a_{1i} \Pf(A_{1i}),
\tag\ebe
$$
where $A_{1i}$ denotes the matrix obtained from $A$ by deleting rows $1$ and $i$, and columns $1$ and $i$. \epf

\endproclaim

{\it Proof of Theorem \tba.} We prove the statement by induction on $k$. For $k=1$ it follows from the fact that $\Pf\left[\matrix 0&a\\-a&0\endmatrix\right]=a$. 

For the induction step, let $k\geq2$ and assume that the statement holds for $k-1$. Let $A$ be the matrix
$$
\align
\left[\matrix
0&\M(G\setminus\{a_1,a_2\})&\M(G\setminus\{a_1,a_3\})&\cdots&\M(G\setminus\{a_1,a_{2k}\})\\
\\
-\M(G\setminus\{a_1,a_2\})&0&\M(G\setminus\{a_2,a_3\})&\cdots&\M(G\setminus\{a_2,a_{2k}\})\\
\\
-\M(G\setminus\{a_1,a_3\})&-\M(G\setminus\{a_2,a_3\})&0&\cdots&\M(G\setminus\{a_3,a_{2k}\})\\
\\
\vdots&\vdots&\vdots& &\vdots\\
\\
-M(G\setminus\{a_1,a_{2k}\})&-\M(G\setminus\{a_2,a_{2k}\})&-\M(G\setminus\{a_3,a_{2k}\})&\cdots&0
\endmatrix\right]\\
\tag\ebf
\endalign
$$
By Lemma {\tbc}, we have
$$
\Pf(A)=\sum_{i=2}^{2k} (-1)^i \M(G\setminus\{a_1,a_i\})\Pf(A_{1i})
\tag\ebg
$$
(recall that $A_{1i}$ is the matrix obtained from $A$ by deleting rows $1$ and $i$, and columns $1$ and $i$). 

Note that the induction hypothesis applied to the graph $G$ and the $2k-2$ vertices in $\overline{\{a_1,a_i\}}$ yields
$$
[\M(G)]^{k-2}\M(G\setminus\overline{\{a_1,a_i\}})=\Pf(A_{1i}),
\tag\ebh
$$
with precisely the same matrices $A_{1i}$ as in (\ebg). It follows from (\ebg) and (\ebh) that
$$
\Pf(A)=[\M(G)]^{k-2}\sum_{i=2}^{2k} (-1)^i \M(G\setminus\{a_1,a_i\})\M(G\setminus\overline{\{a_1,a_i\}}).
\tag\ebi
$$
However, by Proposition {\tbb}, the sum above equals $\M(G)\M(G\setminus\{a_1,\dotsc,a_{2k}\})$. Thus (\ebi) implies (\eba). \epf

\flushpar
\medskip
{\bf Remark 1.} One special situation is when the face $F$ of $G$ containing the vertices $a_1,\dotsc,a_{2k}$ has some pending edges pointing to its interior (see Figure {\fba}), at least one of which has both endpoints in the set $\{a_1,\dotsc,a_{2k}\}$. For defineteness, suppose that $\{a_1,a_2\}$ is such a pending edge, with $a_2$ having degree one. The unusual feature of this situation is that as one moves cyclically around the vertices of $F$, the vertex $a_1$ is visited twice --- once just before encountering $a_2$, and once just after that. Then circular order on $a_1,\dotsc,a_{2k}$ is defined by simply ignoring each such second visit. In fact, one readily checks that if $a_1,\dotsc,a_{2k}$ occur in this circular order {\it modulo the ordering within the endpoints of each pending edge} (i.e., if in the order of the previous statement one is allowed to swap the endpoints of any pending edge connecting two $a_i$'s), the statement of Theorem {\tba} holds without change.

\topinsert
\centerline{\mypic{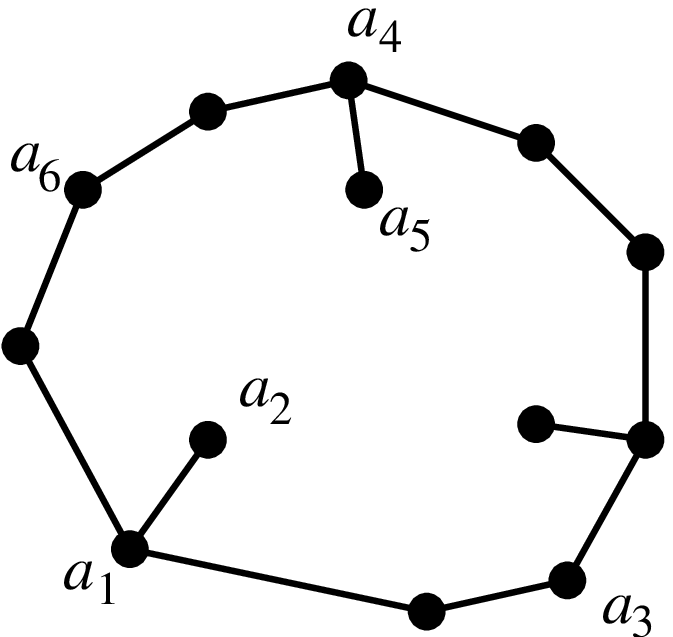}}
\medskip
\centerline{{\smc Figure~{\fba}. {\rm Circular order on a face with pending edges.}}}
\endinsert

\proclaim{Corollary \tbd} Let $G$ be a planar bipartite graph with the same number of vertices in its two color classes. Let the vertices $a_1,\dotsc,a_k,b_k,\dotsc,b_1$ appear in that cyclic order on a face of $G$, and suppose that all the $a_i$'s belong to one color class, and all the $b_j$'s to the other. Then
$$
[\M(G)]^{k-1}\M(G\setminus\{a_1,\dotsc,a_k,b_1,\dotsc,b_k\})
=
\det\left[\M(G\setminus\{a_i,b_j\})\right]_{1\leq i,j\leq k}.
\tag\ebj
$$

\endproclaim

\pf Use Theorem {\tba} to express the left hand side of (\ebj) as the Pfaffian of a $2k\times2k$ matrix. Since $G$ is bipartite, the top left and bottom right quarters of this matrix consist of 0's. Furthermore, if
$$
B=\left[\M(G\setminus\{a_i,b_{k-j+1}\})\right]_{1\leq i,j\leq k}
$$
is the top right quarter, then the bottom left quarter is $-B^T$. Since the Pfaffian is the square root of the determinant, (\ebj) follows, up to sign. The sign turns out to be precisely offset by reversing the order of the columns in the determinant (see e.g. \cite{\Fulmek,Corollary\,1}). This completes the proof. \epf

\mysec{3. A conceptual proof of Eisenk\"olbl's theorem}

In \cite{\Eisen}, Eisenk\"olbl proved the following formula for the number of lozenge tilings of a hexagon with three unit dents along alternating sides.

Recall that the Pochhammer symbol $(a)_k$ is defined by
$$
(a)_k=a(a+1)\cdots(a+k-1).
\tag\eca
$$

\proclaim{Theorem \tca} Let $H_{a,b,c}^{r,s,t}$ be the region obtained from the hexagon of side lengths $a,b+3,c,a+3,b,c+3$ (clockwise, starting with the northern side) by deleting three up-pointing unit triangles from along its boundary as indicated in Figure {\fca}. Then we have
$$
\spreadlines{3\jot}
\align
&
\M(H_{a,b,c}^{r,s,t})=
(r+1)_b(s+1)_c(t+1)_a(a+3-r)_c(b+3-s)_a(c+3-t)_b
\\
&
\times
\frac{\prod_{k=0}^a k!\prod_{k=0}^b k!\prod_{k=0}^c k!\prod_{k=0}^{a+b+c+2} k!}
{\prod_{k=0}^{b+c+2} k!\prod_{k=0}^{a+c+2} k!\prod_{k=0}^{a+b+2} k!}
\\
&
\times
\left[
(a+1)(b+1)(c+1)(a+2-r)(b+2-s)(c+2-t)+(a+1)(b+1)(c+1)rst\right.
\\
&
-(a+2-r)(b+2-s)(c+2-t)rst+(a+1)(c+1)(b+2-s)(c+2-t)rs
\\
&
\left.
+(b+1)(a+1)(a+2-r)(c+2-t)st+(c+1)(b+1)(a+2-r)(b+2-s)rt\right].
\tag\ecb
\endalign
$$

\endproclaim

Our original observation which sparked the current paper was that the 6-term factor above can be written in terms of a 3 by 3 determinant as
$$
\spreadmatrixlines{2\jot}
rstr's't'(a+1)(b+1)(c+1)
\det\left[\matrix
\frac{1}{r'}&\frac{1}{b+1}&\frac{1}{t}\\
-\frac{1}{r}&\frac{1}{s'}&\frac{1}{c+1}\\
-\frac{1}{a+1}&-\frac{1}{s}&\frac{1}{t'}\\
\endmatrix\right],
$$
where for brevity of notation we wrote $r'=a+2-r$, $s'=b+2-s$ and $t'=c+2-t$.

The reason for this is apparent from our proof of Eisenk\"olbl's theorem, presented below.

\topinsert
\centerline{\mypic{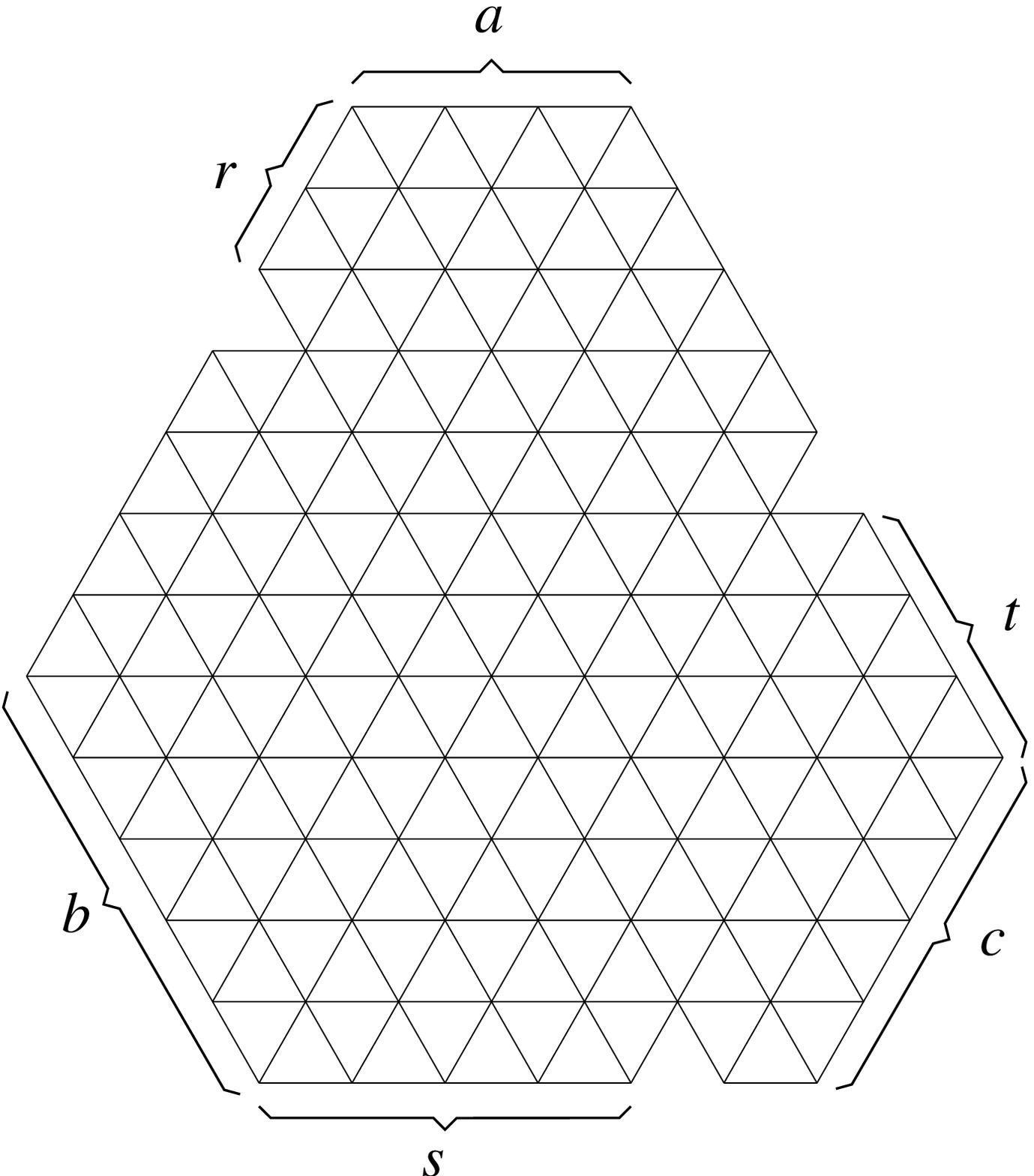}}
\centerline{Figure~{\fca}. {\rm $H_{3,4,5}^{4,3,2}$.}}
\endinsert

\medskip
\pf Denote by $H^*_{a,b,c}$ the region obtained from $H_{a,b,c}^{r,s,t}$ by filling back the three unit dents along its sides, and adding three additional unit triangles sticking out next to the bottom left, right, and top left corners as indicated in Figure {\fcb}. 

Apply Theorem {\tba} to the planar dual graph\footnote{ The planar dual graph of a region $R$ on the triangular lattice is the graph whose vertices are the unit triangles in $R$, and whose edges connect vertices corresponding to unit triangles that share an edge.} $G$ of $H^*_{a,b,c}$, with $k=3$, and the six removed vertices chosen to correspond to the three dents in the statement of the theorem and the three unit triangles that stick out \footnote{ Note that, by Remark 1, we do not need to treat separately the cases when some of $r$, $s$, $t$ are 0; if that happens, we have a ``pending edge'' situation, and Theorem {\tba} still applies.}(see Figure {\fcc}).  Let $a_1$, $a_2$ and $a_3$ be the dents along the sides of lengths $a+3$, $b+3$ and $c+3$, respectively, and let $b_1$, $b_2$ and $b_3$ be the unit triangles that stick out from the corresponding edges. Then $b_1,a_1,b_2,a_2,b_3,a_3$ occur in cyclic order along the unbounded face of $G$ (see the convention in Remark 1 for the case when an $a_i$ shares an edge with a $b_j$). Thus, by Theorem {\tba}, we obtain that $[\M(G)]^2\M(G\setminus\{b_1,a_1,b_2,a_2,b_3,a_3\}$ is equal to the Pfaffian of the matrix
$$
\spreadmatrixlines{2\jot}
\left[\matrix
0&\M(G_{b_1,a_1})&\M(G_{b_1,b_2})&\M(G_{b_1,a_2})&\M(G_{b_1,b_3})&\M(G_{b_1,a_3})\\
-\M(G_{b_1,a_1})&0&\M(G_{a_1,b_2})&\M(G_{a_1,a_2})&\M(G_{a_1,b_3})&\M(G_{a_1,a_3})\\
-\M(G_{b_1,b_2})&-\M(G_{a_1,b_2})&0&\M(G_{b_2,a_2})&\M(G_{b_2,b_3})&\M(G_{b_2,a_3})\\
-\M(G_{b_1,a_2})&-\M(G_{a_1,a_2})&-\M(G_{b_2,a_2})&0&\M(G_{a_2,b_3})&\M(G_{a_2,a_3})\\
-\M(G_{b_1,b_3})&-\M(G_{a_1,b_3})&-\M(G_{b_2,b_3})&-\M(G_{a_2,b_3})&0&\M(G_{b_3,a_3})\\
-\M(G_{b_1,a_3})&-\M(G_{a_1,a_3})&-\M(G_{b_2,a_3})&-\M(G_{a_2,a_3})&-\M(G_{b_3,a_3})&0\\
\endmatrix\right],
$$

\topinsert
\twoline{\mypic{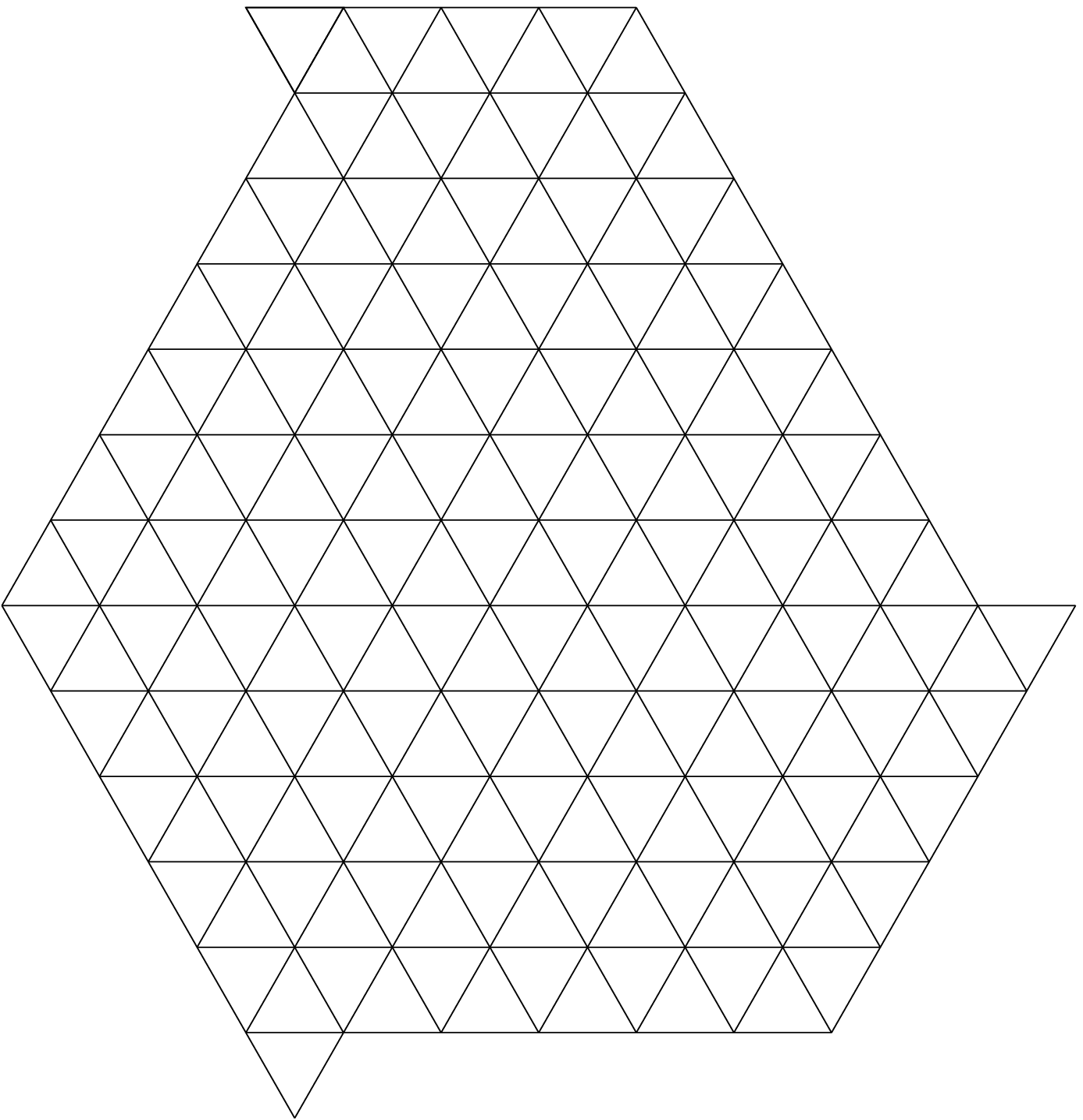}}{\mypic{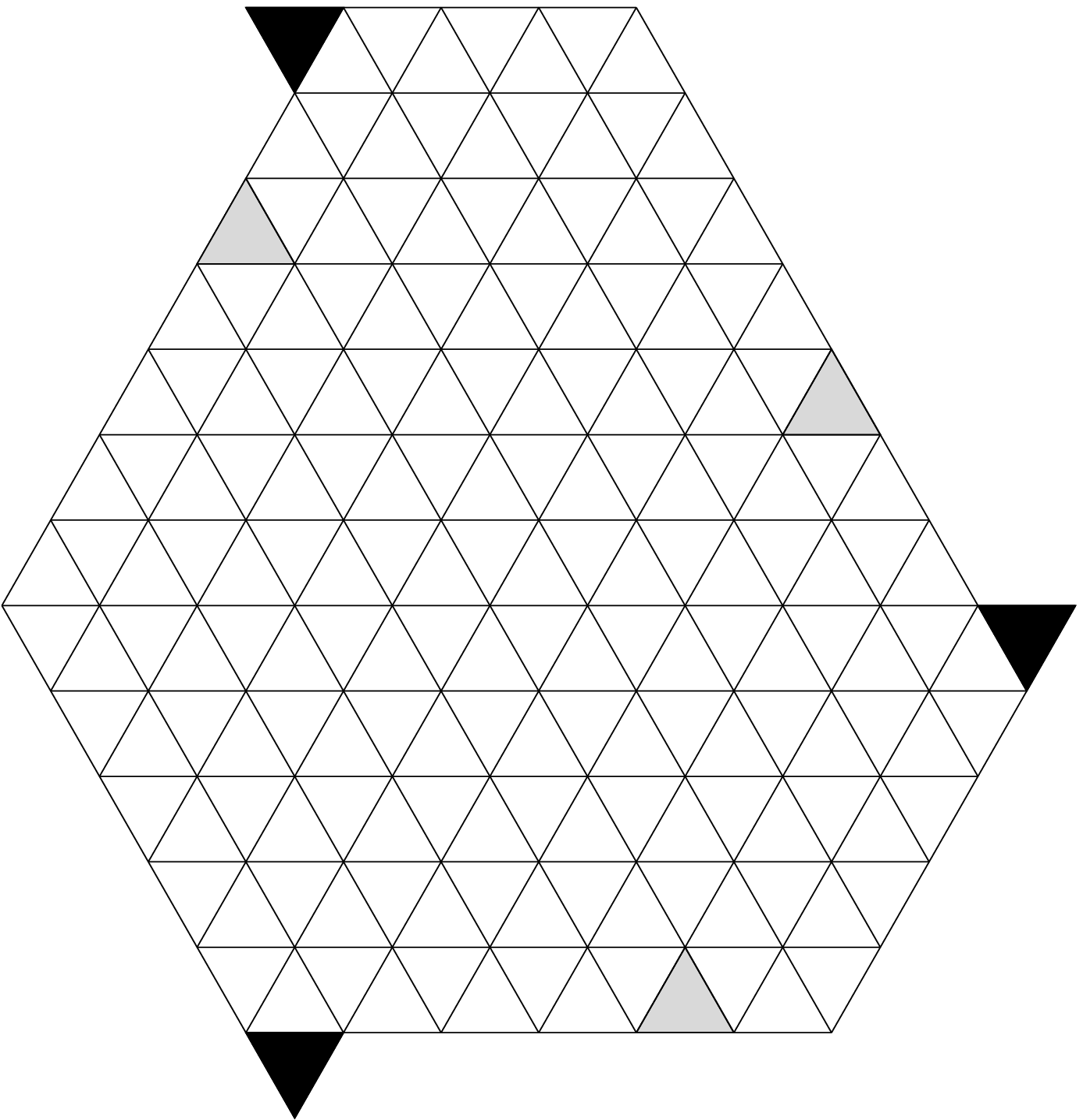}}
\twoline{Figure~{\fcb}. {\rm $H_{3,4,5}^*$.}}{Figure~{\fcc}. {\rm Choosing the vertices.}}
\endinsert

\flushpar
where for brevity of notation we wrote $G_{u,v}$ for $G\setminus\{u,v\}$. However, since $G$ is bipartite with the same number of vertices in the two color classes, all entries in the above matrix corresponding to removing two $a_i$'s or two $b_j$'s are zero. Thus we obtain
$$
\spreadlines{3\jot}
\spreadmatrixlines{2\jot}
\align
&
[\M(G)]^2\M(G\setminus\{b_1,a_1,b_2,a_2,b_3,a_3\})=
\\
&
\Pf
\left[\matrix
0&\M(G_{b_1,a_1})&0&\M(G_{b_1,a_2})&0&\M(G_{b_1,a_3})\\
-\M(G_{b_1,a_1})&0&\M(G_{a_1,b_2})&0&\M(G_{a_1,b_3})&0\\
0&-\M(G_{a_1,b_2})&0&\M(G_{b_2,a_2})&0&\M(G_{b_2,a_3})\\
-\M(G_{b_1,a_2})&0&-\M(G_{b_2,a_2})&0&\M(G_{a_2,b_3})&0\\
0&-\M(G_{a_1,b_3})&0&-\M(G_{a_2,b_3})&0&\M(G_{b_3,a_3})\\
-\M(G_{b_1,a_3})&0&-\M(G_{b_2,a_3})&0&-\M(G_{b_3,a_3})&0\\
\endmatrix\right]
\\
\tag\ecc
\endalign
$$
Reordering rows and columns --- each simultaneous interchange of two rows and the corresponding two columns results in a sign change for the Pfaffian --- we obtain from (\ecc) that
$$
\spreadmatrixlines{2\jot}
[\M(G)]^2\M(G\setminus\{b_1,a_1,b_2,a_2,b_3,a_3\})=
-\Pf
\left[\matrix
0&B\\
-B^T&0\\
\endmatrix\right],
\tag\ecdd
$$
where
$$
\spreadmatrixlines{2\jot}
B=
\left[\matrix
\M(G\setminus\{b_1,a_1\})&\M(G\setminus\{b_1,a_2\})&\M(G\setminus\{b_1,a_3\})\\
-\M(G\setminus\{b_2,a_1\})&\M(G\setminus\{b_2,a_2\})&\M(G\setminus\{b_2,a_3\})\\
-\M(G\setminus\{b_3,a_1\})&-\M(G\setminus\{b_3,a_2\})&\M(G\setminus\{b_3,a_3\})\\
\endmatrix\right].
$$
Since for any $k\times k$ matrix $C$
$$
\Pf\left[\matrix 0&C\\-C^T&0\endmatrix\right]=(-1)^{k(k-1)/2}\det(C)
\tag\ecddd
$$
(see e.g. \cite{\Fulmek,Corollary\,1}), we obtain from (\ecdd) that
$$
\spreadlines{3\jot}
\spreadmatrixlines{2\jot}
\align
&\!\!\!\!\!\!\!\!\!\!\!\!
[\M(G)]^2\M(G\setminus\{b_1,a_1,b_2,a_2,b_3,a_3\})=
\\
&\ \ \ \ \ \ \ \ \ \ \ \ 
\det
\left[\matrix
\M(G\setminus\{b_1,a_1\})&\M(G\setminus\{b_1,a_2\})&\M(G\setminus\{b_1,a_3\})\\
-\M(G\setminus\{b_2,a_1\})&\M(G\setminus\{b_2,a_2\})&\M(G\setminus\{b_2,a_3\})\\
-\M(G\setminus\{b_3,a_1\})&-\M(G\setminus\{b_3,a_2\})&\M(G\setminus\{b_3,a_3\})\\
\endmatrix\right].
\\
\tag\ecd
\endalign
$$
Note that by the way we set things up, $\M(G\setminus\{b_1,a_1,b_2,a_2,b_3,a_3\})$ is precisely the left hand side of (\ecb), which we want to determine. The fortunate situation is that all the remaining perfect matching counts in (\ecd) can be readily obtained.

\topinsert
\centerline{\mypic{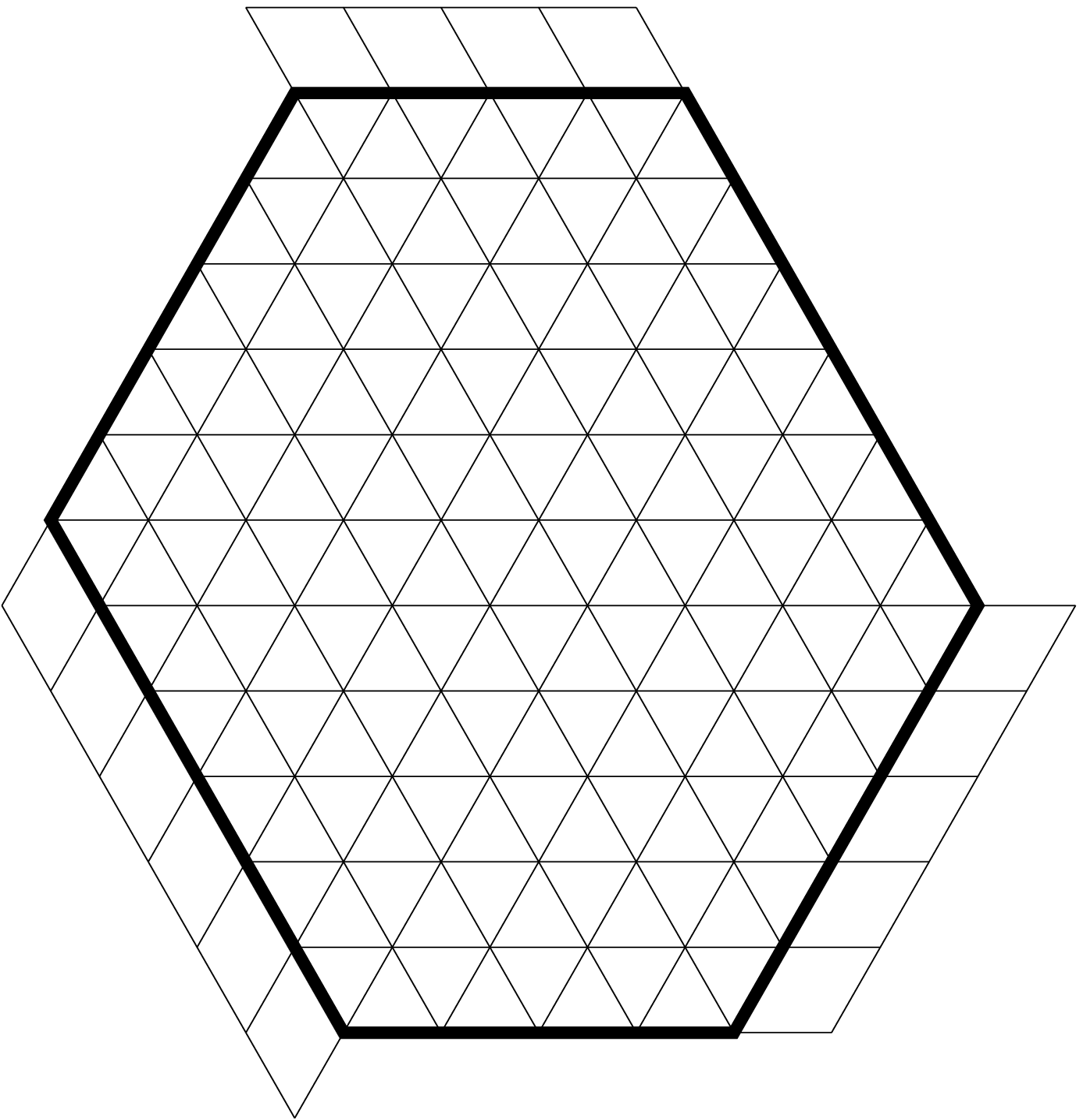}}
\medskip
\centerline{{\smc Figure~{\fcd}. {\rm Removing the forced tiles in $H^*_{a,b,c}$.}}}
\endinsert

\topinsert
\twoline{\mypic{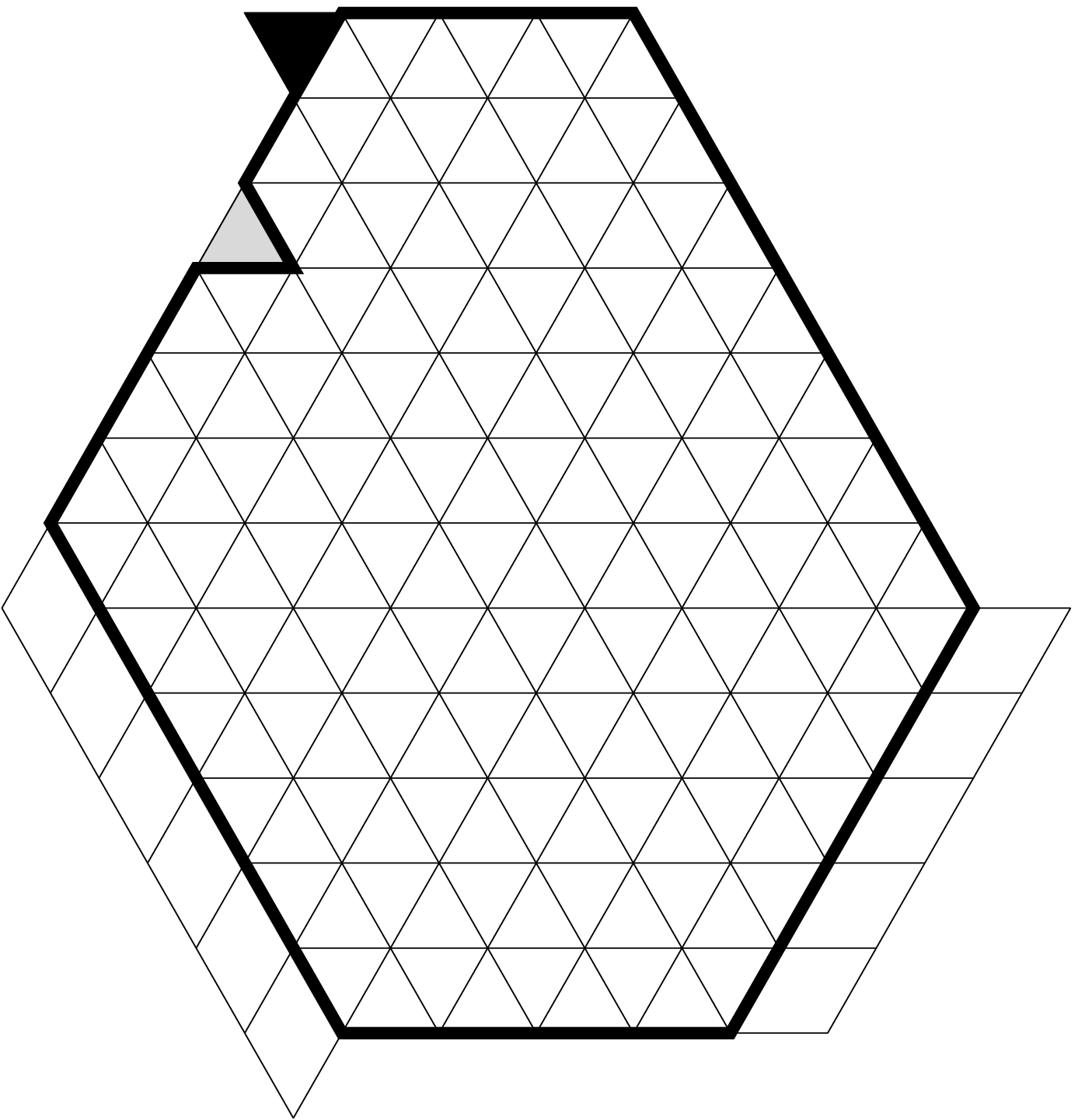}}{\mypic{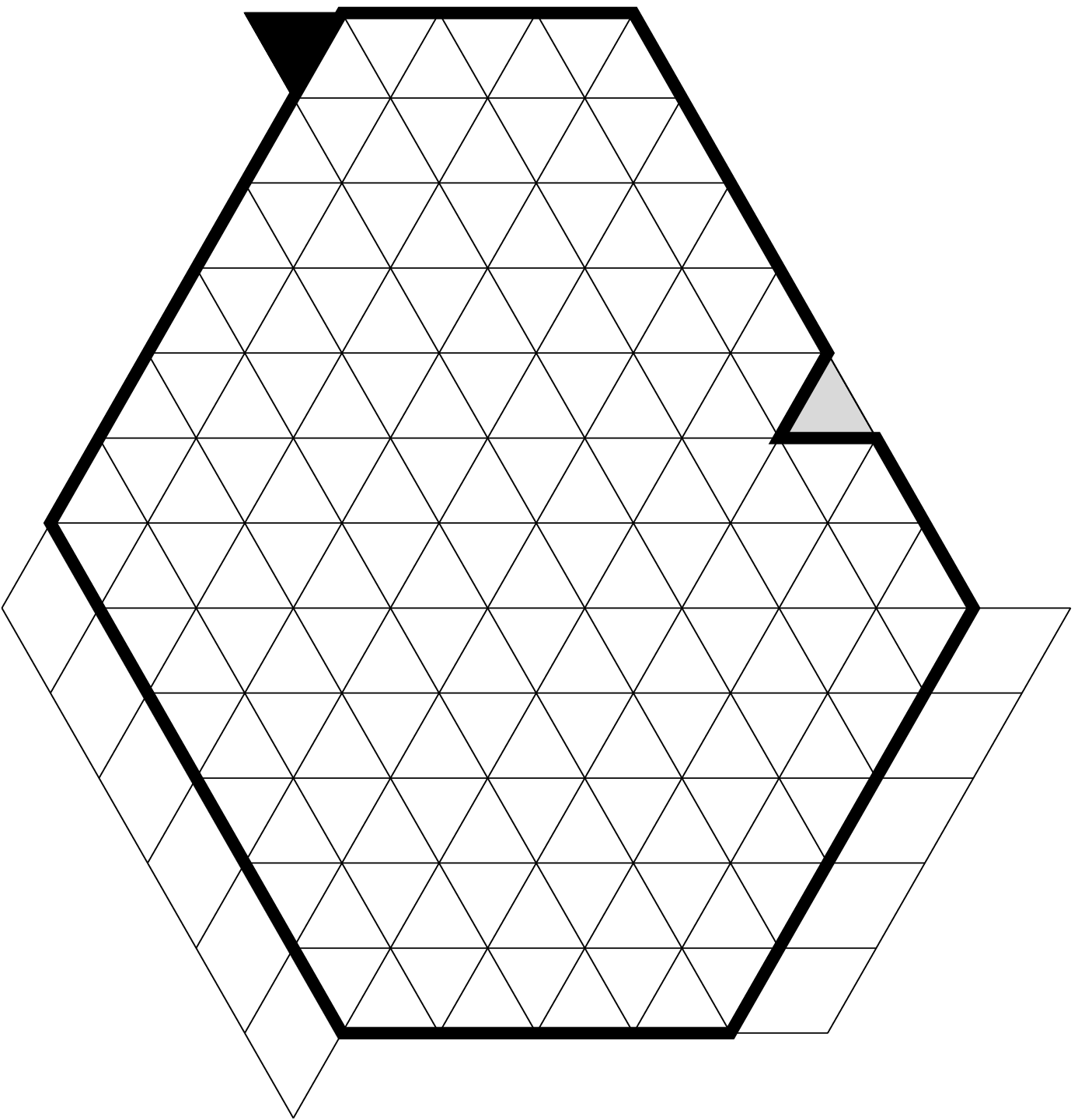}}
\medskip
\centerline{Figure~{\fce}. {\rm Removing forced tiles in the $H_{a,b,c}^*\setminus\{a_i,b_j\}$'s.}}
\endinsert

Indeed, due to forced tiles, we have that
$$
\M(G)=\M(H^*_{a,b,c})=\M(H_{a+1,b+1,c+1}),
$$
where $H_{a+1,b+1,c+1}$ is the hexagon of sides $a+1,b+1,c+1,a+1,b+1,c+1$ (clockwise from top; see Figure {\fcd}); hence by MacMahon's classical theorem on boxed plane partitions \cite{\MacM} (which are well-known to be equivalent to lozenge tilings of hexagons) we have
$$
\M(H_{a+1,b+1,c+1})=
\frac{\prod_{k=0}^a k!\prod_{k=0}^b k!\prod_{k=0}^c k!\prod_{k=0}^{a+b+c+2} k!}
{\prod_{k=0}^{b+c+1} k!\prod_{k=0}^{a+c+1} k!\prod_{k=0}^{a+b+1} k!}.
\tag\ece
$$
Furthermore, for all $i,j\in\{1,2,3\}$, after removing the forced tiles, the region corresponding to $G\setminus\{a_i,b_j\}$ is a hexagon with a single unit dent along one of its sides (Figure {\fce} illustrates the two types of regions that arise this way). The number of its lozenge tilings follows thus from the general formula in Lemma {\tcb} below --- in addition to an isolated dent around the middle, include consecutive runs of $x_i$'s at the extreme left and right in Figure {\fcf} to turn the trapezoidal region into a hexagon with a single dent (see Figure {\fcg} for an illustration of this).
It is routine to verify that plugging in the above explicit formulas into (\ecd) one obtains~(\ecb).~\epf

\topinsert
\twoline{\mypic{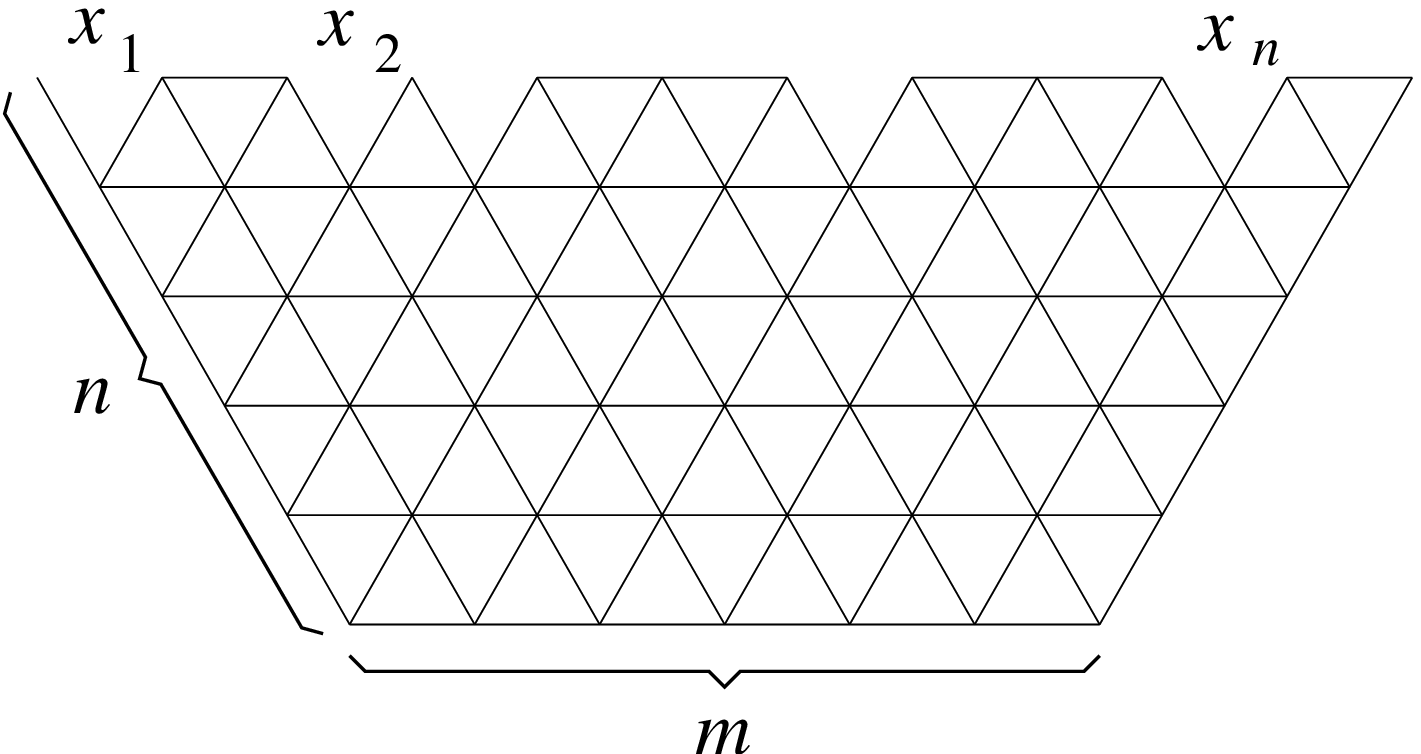}}{\mypic{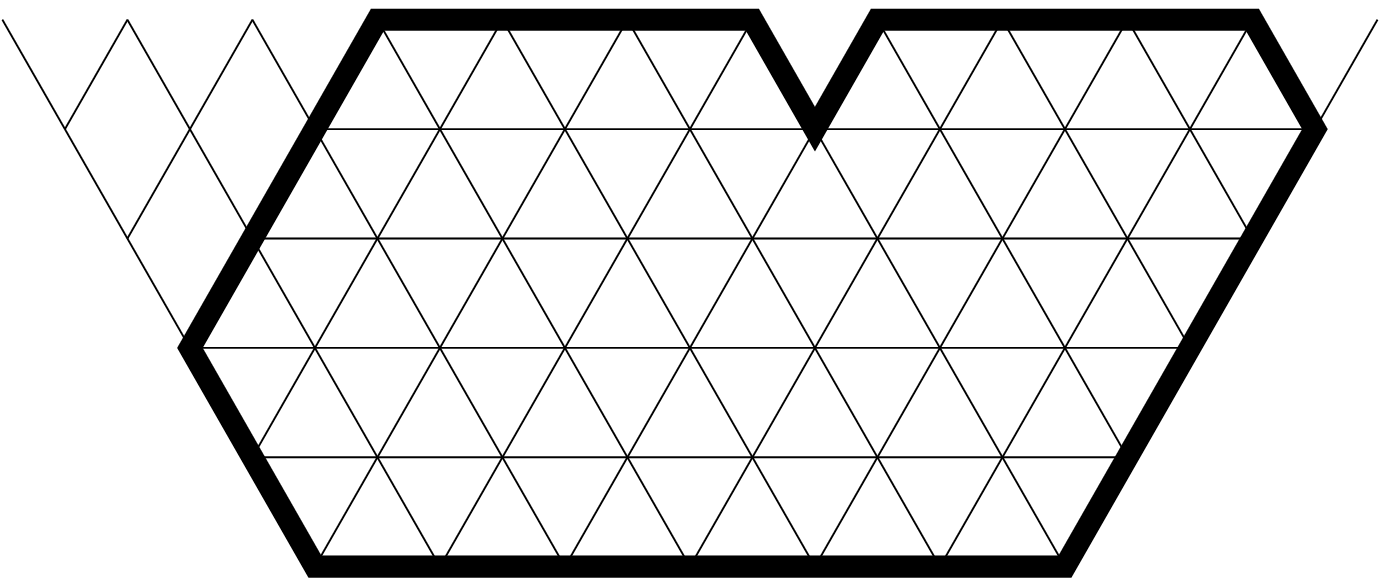}}
\twoline{Figure~{\fcf}. {\rm $T_{6,5}(1,3,4,7,10)$.}}{Figure~{\fcg}. {\rm Single dent hexagon as a $T$-region.}}

\endinsert

The following formula is Cohn, Larsen and Propp's \cite{\CLP} translation to lozenge tilings of a classical result of Gelfand and Tsetlin \cite{\GT}.

\proclaim{Proposition \tcb} Let $T_{m,n}(x_1,\dotsc,x_n)$ be the region obtained from the trapezoid of side lengths $m$, $n$, $m+n$, $n$ (clockwise from bottom) by removing the down-pointing unit triangles from along its top that are in positions $x_1,x_2,\dotsc,x_n$ as counted from left to right. Then
$$
\M(T_{m,n}(x_1,\dotsc,x_n))=\prod_{1\leq i<j\leq n}\frac{x_j-x_i}{j-i}.
\tag\ecf
$$

\endproclaim

\mysec{4. A generalization}

We generalize Eisenk\"olbl's regions $H_{x,y,z}^{r,s,t}$ of the previous section as follows. 
Let $H_{x,y,z}^k$ be the hexagon on the triangular lattice whose sides have lengths $x,y+k,z,x+k,y,z+k$, in clockwise order starting at the top. There are precisely $x+y+z+3k$ up-pointing unit lattice triangles in it that share an edge with the boundary --- $x+k$, $y+k$, {resp.} $z+k$ along the southern, northeastern, {resp.} northwestern sides. Choose $k$ of them, and denote them by $a_1,\dotsc,a_k$. 
Our generalization of Eisenk\"olbl's regions is the family of regions of type $H_{x,y,z}^k\setminus\{a_1,\dotsc,a_k\}$ (see Figure {\fda} for an example).

\proclaim{Theorem \tda} Let $H_{x,y,z}^\star$ be the region obtained from $H_{x,y,z}^k$ by augmenting it with three strings of contiguous down-pointing unit triangles along its boundary as shown in Figure~${\fdb}({\text{\rm a}})$; the length of the string on each side is equal to the number of $a_i$'s in $H_{x,y,z}^k\setminus\{a_1,\dotsc,a_k\}$ along that side. Denote the $k$ down-pointing unit triangles in these strings by $b_1,\dotsc,b_k$.
Let $c_1,\dotsc,c_{2k}$ be the elements of the set $\{a_1,\dotsc,a_k\}\cup\{b_1,\dotsc,b_k\}$ listed in a cyclic order\footnote{ If $a_1$ (resp., $b_1$) is the leftmost $a_i$ (resp., $b_i$) along the bottom side in Figure {\fdb}(b), and $a_1,\dotsc,a_7$ (resp., $b_1,\dotsc,b_7$) occur in counterclockwise order, then one such cyclic order of the union of the $a_i$'s and $b_i$'s is for instance $b_1,b_2,a_1,b_4,a_2,a_3,a_4,b_5,a_5,b_6,b_7,a_6,a_7$.}, as explained in Remark 1.
Then we have
$$
\M(H_{x,y,z}^k\setminus\{a_1,\dotsc,a_k\}) 
=
\frac{1}{\left[\M(H_{x,y,z}^\star)\right]^{k-1}}
\Pf\left[
\M(H_{x,y,z}^\star\setminus\{c_i,c_j\})
\right]_{1\leq i<j\leq 2k}
\tag\eda
$$
where the quantities on the right hand side are given by explicit formulas: $\M(H_{x,y,z}^\star)$ by $(\ece)$, $\M(H_{x,y,z}^\star\setminus\{a_i,b_j\})$ by $(\ecf)$ if $a_i$ and $b_j$ are along the same side and by Proposition~${\tdb}$  if $a_i$ and $b_j$ are along different sides, and $\M(H_{x,y,z}^\star\setminus\{a_i,a_j\})=\M(H_{x,y,z}^\star\setminus\{b_i,b_j\})=0$.

\endproclaim

\topinsert
\centerline{\mypic{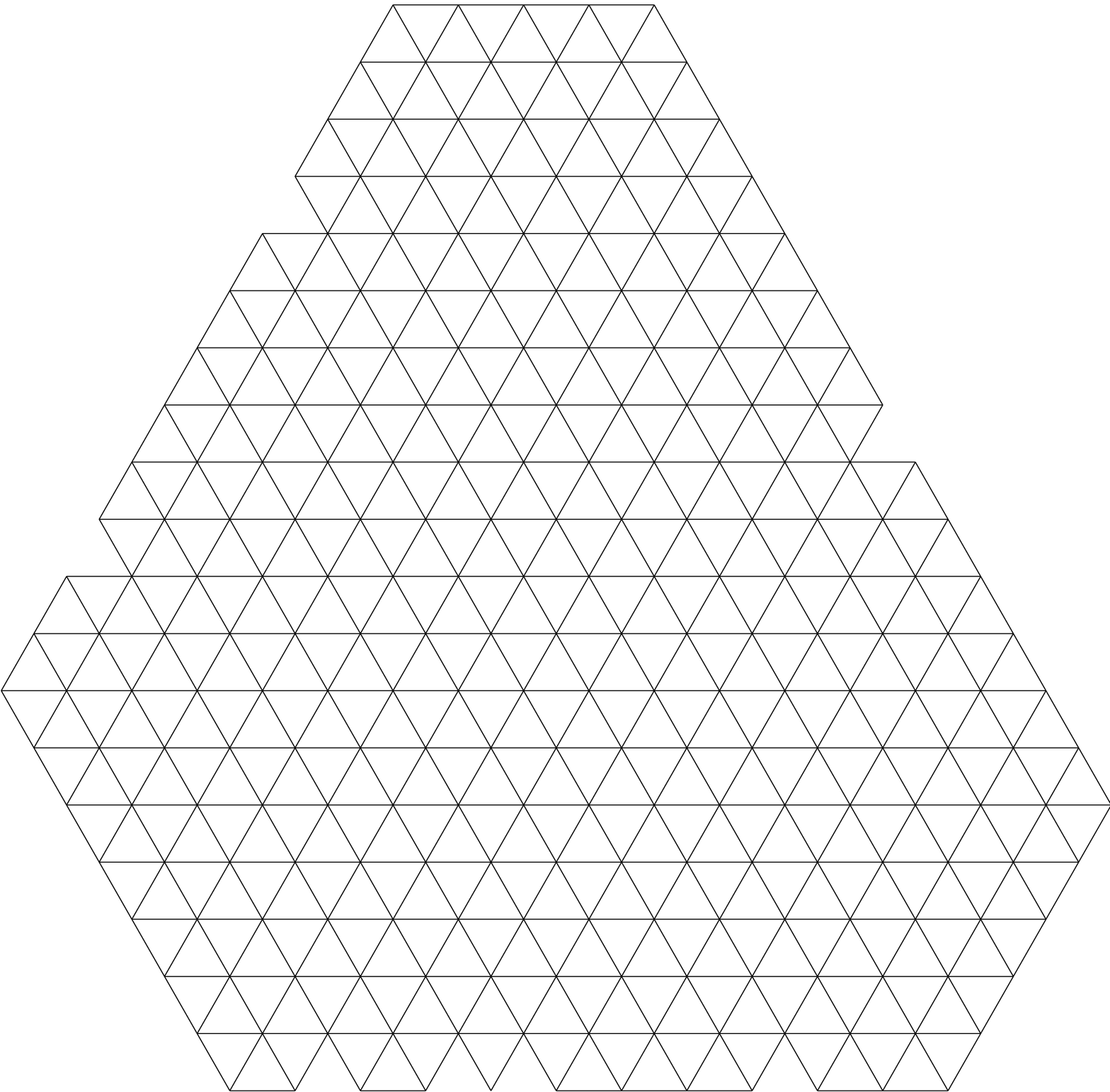}}
\medskip
\centerline{{\smc Figure~{\fda}. {\rm A hexagon with seven up-pointing dents.}}}
\endinsert

\topinsert
\twolinetwo{\mypic{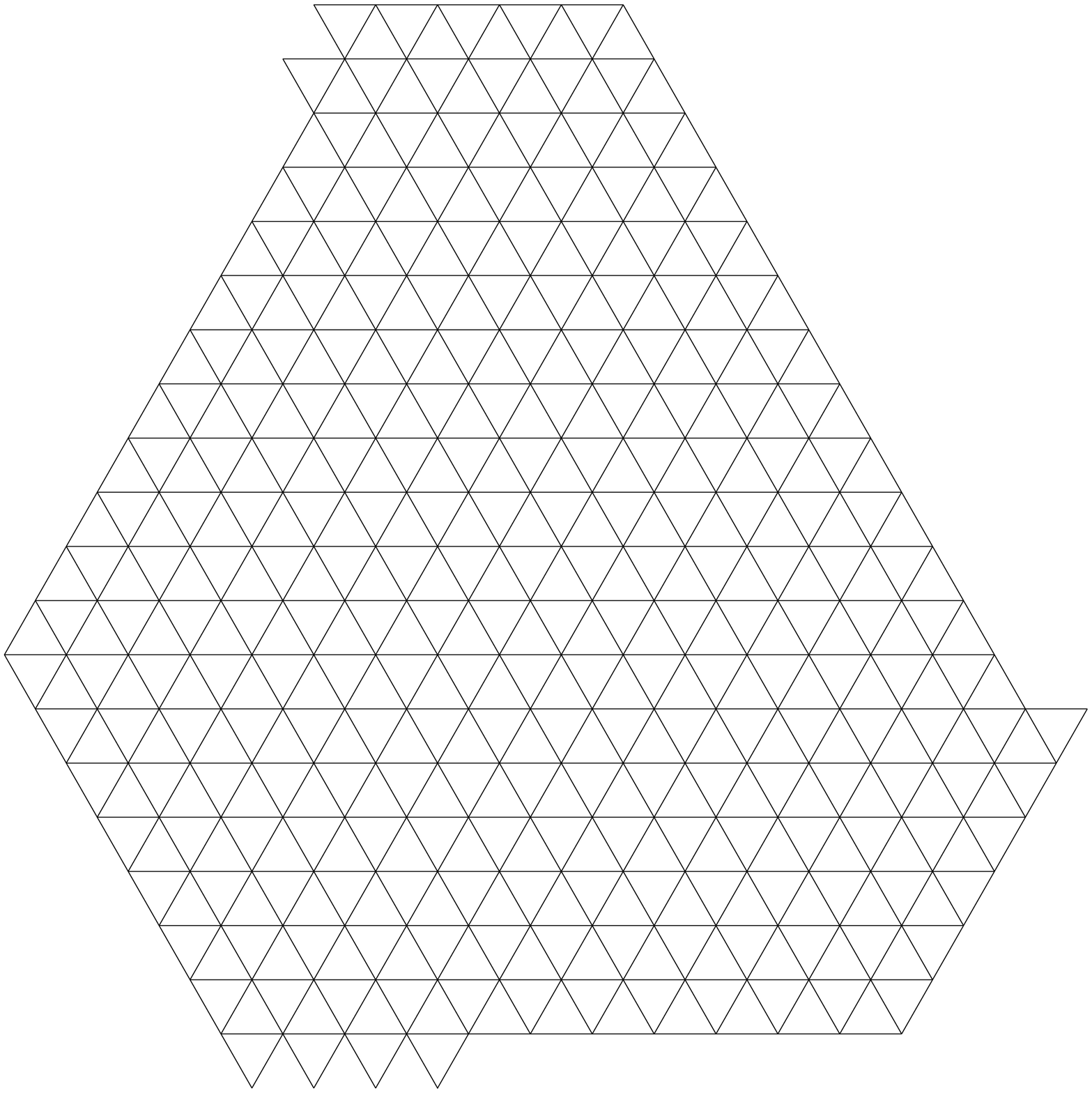}}{\mypic{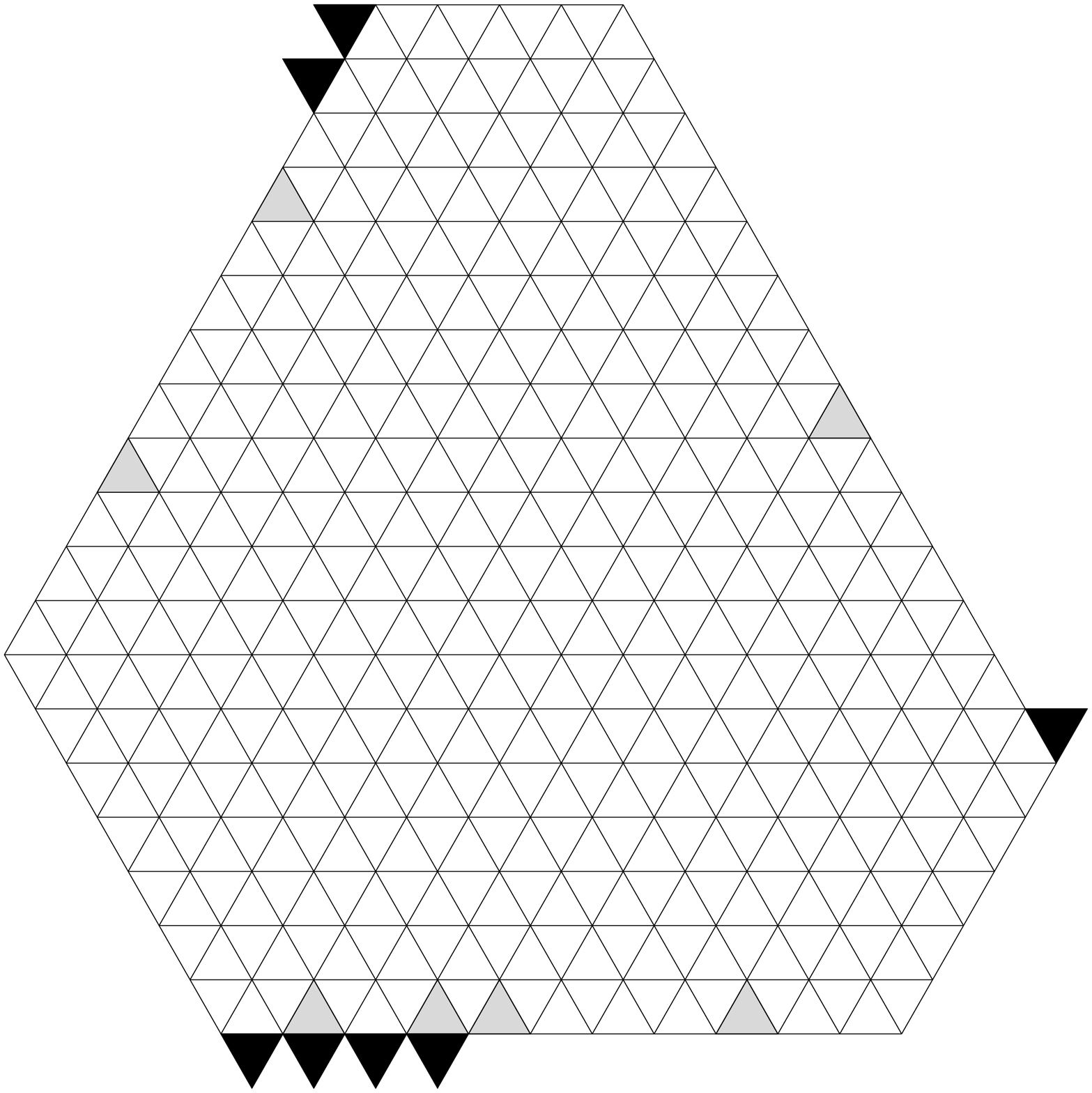}}
\twoline{\rm (a)\ \ \ \ \ \ \ \ \ }{\rm \ \ \ \ \ \ \ \ \ \ \ \ \ \ \ \ \ \ \ \ \ \ \ (b)}
\medskip
\centerline{Figure~{\fdb}. \!{\rm (a). Region to which we apply condensation. (b).  Choosing the vertices.}}
\endinsert

{\bf Remark 2.} After simultaneous reorderings of rows and columns (which preserve the Pfaffian, up to sign), the matrix in (\eda) can always be brought to the form
$$
\left[\matrix 0&B\\-B^T&0\endmatrix\right],
$$
where $B$ is a $k\times k$ matrix. Thus, by (\ecddd),
formula (\eda) yields a determinant expression for $\M(H_{x,y,z}^k\setminus\{a_1,\dotsc,a_k\})$.  

The subtlety is that the entries of $B$ are {\it signed} $\M(H_{x,y,z}^\star\setminus\{a_i,b_j\})$'s. For instance, for the example in Figure {\fdb}, the resulting matrix $B$ is
$$
\!\!\!\!\!\!\!\!\!\!\!\!\!\!\!\!
\spreadmatrixlines{2\jot}
B=
\left[\matrix
-m_{a_1,b_1}&-m_{a_1,b_2}&m_{a_1,b_3}&m_{a_1,b_4}&m_{a_1,b_5}&m_{a_1,b_6}&m_{a_1,b_7}\\
-m_{a_2,b_1}&-m_{a_2,b_2}&-m_{a_2,b_3}&-m_{a_2,b_4}&m_{a_2,b_5}&m_{a_2,b_6}&m_{a_2,b_7}\\
-m_{a_3,b_1}&-m_{a_3,b_2}&-m_{a_3,b_3}&-m_{a_3,b_4}&m_{a_3,b_5}&m_{a_3,b_6}&m_{a_3,b_7}\\
-m_{a_4,b_1}&-m_{a_4,b_2}&-m_{a_4,b_3}&-m_{a_4,b_4}&m_{a_4,b_5}&m_{a_4,b_6}&m_{a_4,b_7}\\
-m_{a_5,b_1}&-m_{a_5,b_2}&-m_{a_5,b_3}&-m_{a_5,b_4}&-m_{a_5,b_5}&m_{a_5,b_6}&m_{a_5,b_7}\\
-m_{a_6,b_1}&-m_{a_6,b_2}&-m_{a_6,b_3}&-m_{a_6,b_4}&-m_{a_6,b_5}&-m_{a_6,b_6}&-m_{a_6,b_7}\\
-m_{a_7,b_1}&-m_{a_7,b_2}&-m_{a_7,b_3}&-m_{a_7,b_4}&-m_{a_7,b_5}&-m_{a_7,b_6}&-m_{a_7,b_7}\\
\endmatrix\right],
$$
where for brevity of notation we wrote $m_{a_i,b_j}$ for $\M(H_{x,y,z}^\star\setminus\{a_i,b_j\})$.


\pf Formula (\eda) follows directly from Theorem {\tba}, with $G$ chosen to be the planar dual graph of the region $H_{x,y,z}^\star$, and $a_1,\dotsc,a_k$ and $b_1,\dotsc,b_k$ chosen to be the vertices of $G$ corresponding to the unit triangles $a_1,\dotsc,a_k,b_1,\dotsc,b_k$ in the statement of Theorem~{\tda}.

\topinsert
\twolinetwo{\mypic{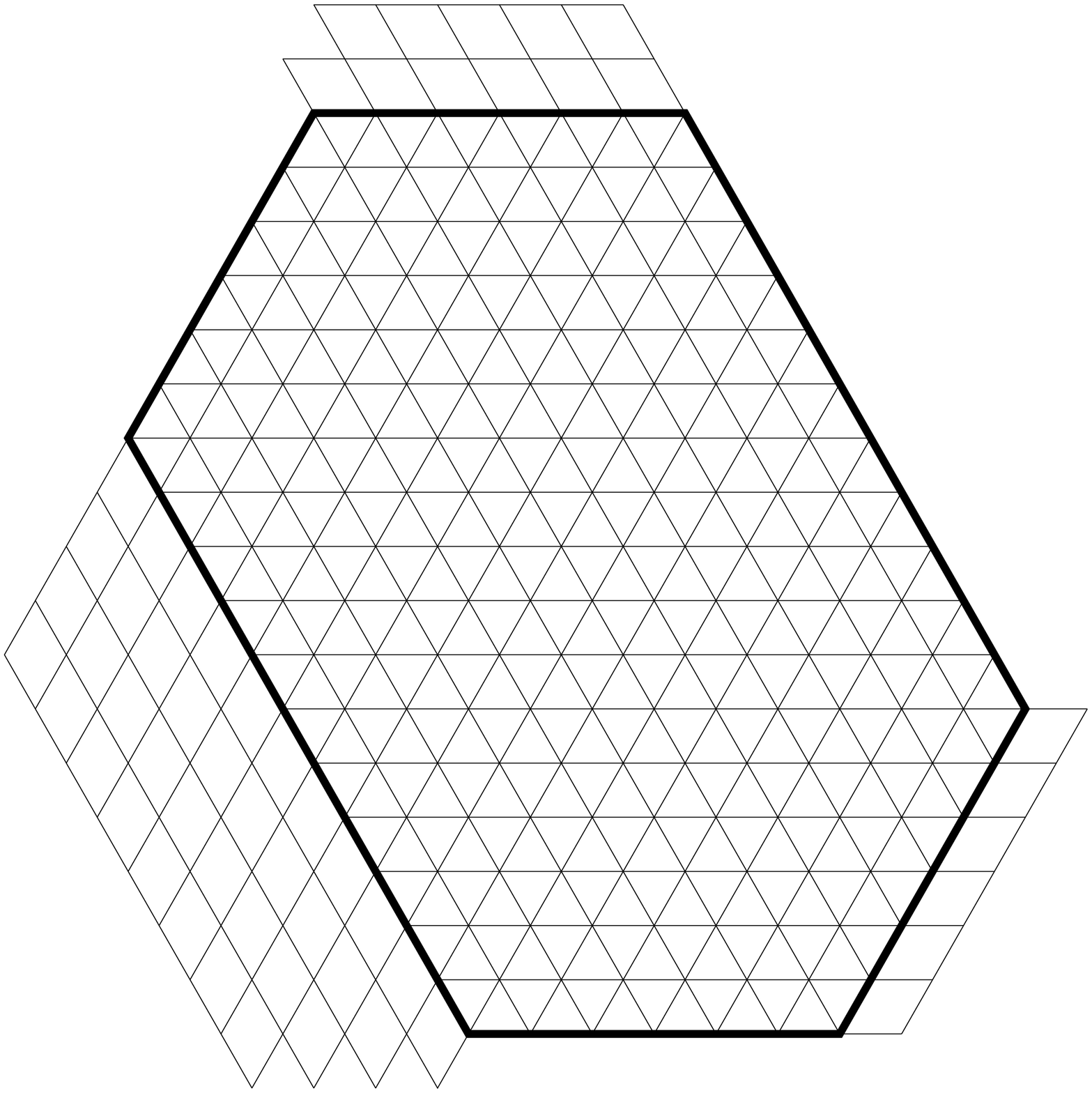}}{\mypic{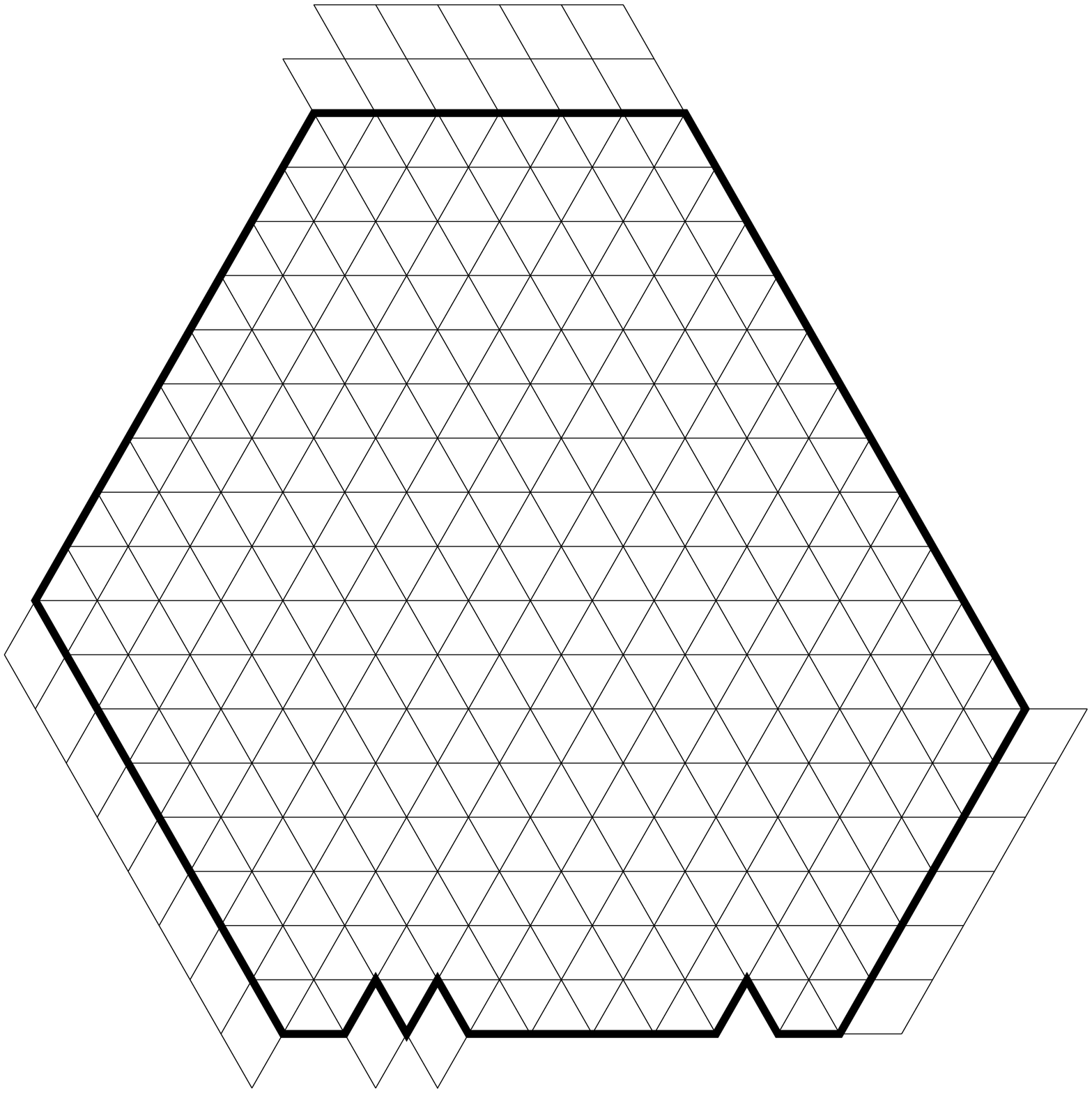}}
\twoline{\rm (a)\ \ \ \ \ \ \ \ \ }{\rm \ \ \ \ \ \ \ \ \ \ \ \ \ \ \ \ \ \ \ \ \ \ \ (b)}
\medskip
\centerline{Figure~{\fdc}. \!{\rm Region obtained by removing the forced lozenges (a) from $H_{x,y,z}^\star$,}}
\centerline{and (b) from $H_{x,y,z}^\star\setminus\{a_i,b_j\}$ when $a_i$, $b_j$ are on the same side.}
\endinsert

\topinsert
\twolinetwo{\mypic{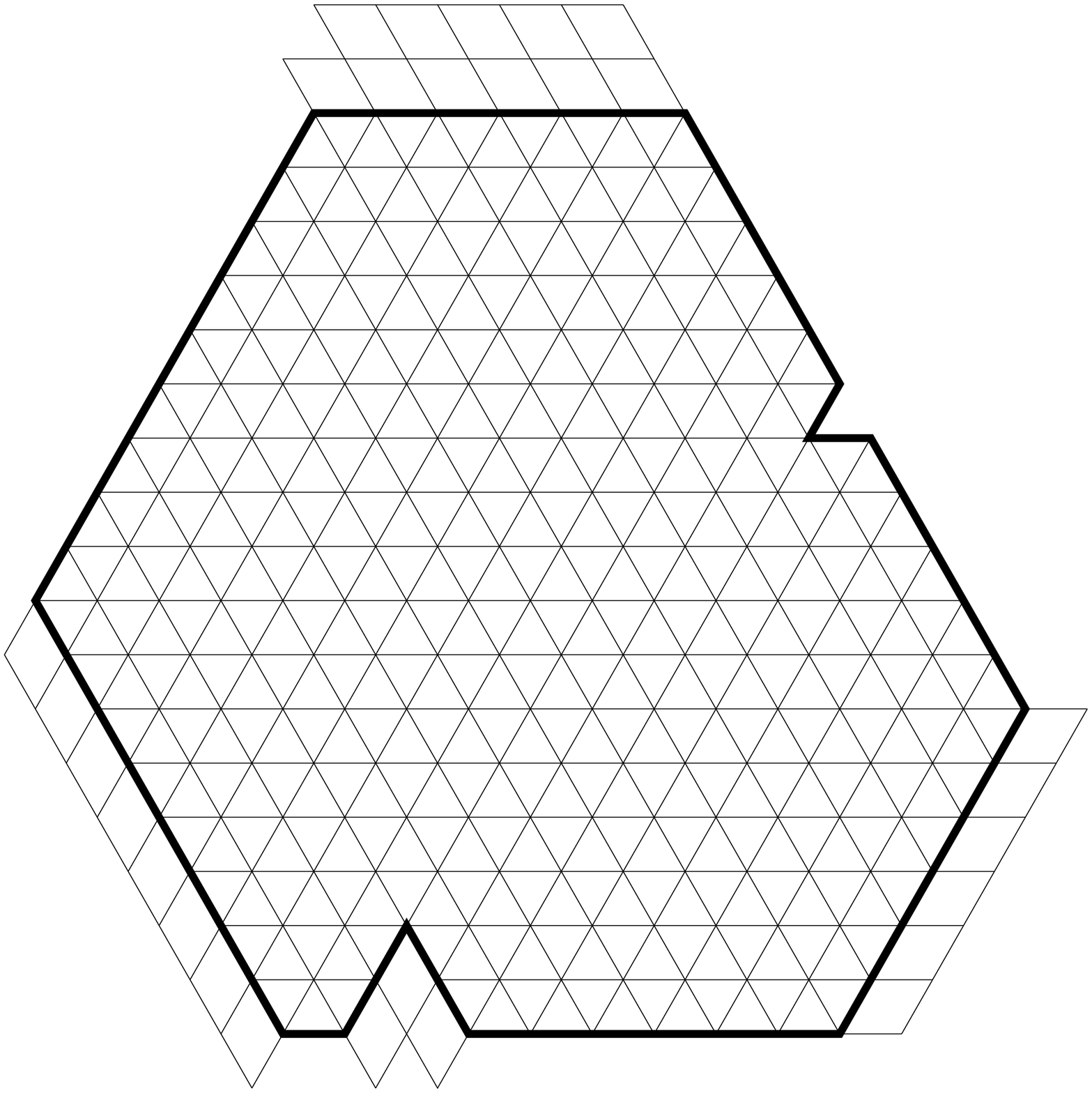}}{\mypic{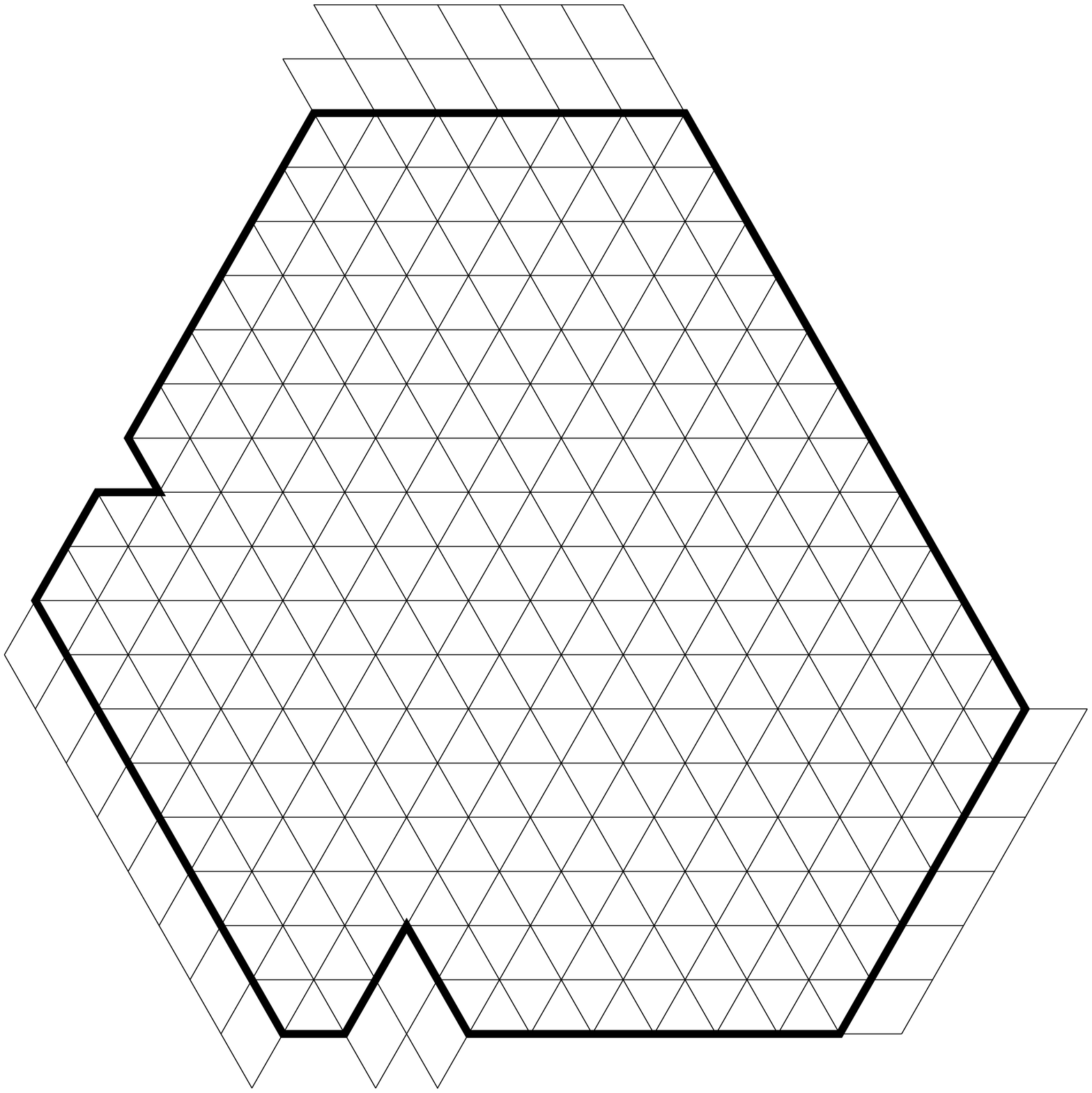}}
\twoline{\rm (a)\ \ \ \ \ \ \ \ \ }{\rm \ \ \ \ \ \ \ \ \ \ \ \ \ \ \ \ \ \ \ \ \ \ \ (b)}
\medskip
\centerline{Figure~{\fdd}. \!{\rm The two types of regions obtained from $H_{x,y,z}^\star\setminus\{a_i,b_j\}$}}
\centerline{by removing the forced lozenges, when $a_i$, $b_j$ are on different sides.}
\endinsert

In order to see how the quantities on the right hand side of (\eda) are given by the indicated formulas, let us consider first the region $H_{x,y,z}^\star$. The three strings of $b_i$'s in it force many lozenges to be part of every tiling of $H_{x,y,z}^\star$. After all these forced lozenges are removed, the resulting region is a lattice hexagon (see Figure {\fdc}(a)). Since $H_{x,y,z}^\star$ is balanced (i.e., has the same number of up-pointing and down-pointing unit triangles), so is the resulting hexagon. Then the lengths of opposite sides must be the same, and the number of lozenge tilings is indeed given by formula (\ece).

We turn next to the entries of type $\M(H_{x,y,z}^\star\setminus\{a_i,b_j\})$, where $a_i$ and $b_j$ are along the same side of $H_{x,y,z}^\star$. Here we distinguish two cases. If in $H_{x,y,z}^\star\setminus\{a_i,b_j\}$, below the removed unit triangle $a_i$, there is an unremoved unit triangle $b_k$, then there is no way to cover $b_k$ by a lozenge, so $\M(H_{x,y,z}^\star\setminus\{a_i,b_j\})=0$ in this case. 

Otherwise, either $a_i$ and $b_j$ share an edge, or $a_i$ does not share an edge with any of the $b_k$'s. Figure {\fdc}(b) illustrates the latter situation. Clearly, after removing the forced lozenges, the resulting region is of the type covered by Proposition {\tcb}, so $\M(H_{x,y,z}^\star\setminus\{a_i,b_j\})$ is given in this case by formula (\ecf). One readily sees that the same holds in the former situation.

The remaining entries of type $\M(H_{x,y,z}^\star\setminus\{a_i,b_j\})$ are those for which $a_i$ and $b_j$ were removed from along different sides of $H_{x,y,z}^\star$. There are two different situations to distinguish, corresponding to the cases when the side from which $a_i$ was removed is the next nearest neighbor of the side from which $b_j$ was removed in the counter-clockwise direction, or in the clockwise direction (these are illustrated in Figures {\fdd}(a) and (b), respectively). After removing the forced lozenges, the resulting regions are readily seen to be of the types covered by Proposition {\tdb}(a) and (b), respectively.

Since $H_{x,y,z}^\star\setminus\{a_i,a_j\}$ and $H_{x,y,z}^\star\setminus\{b_i,b_j\}$ are not balanced (i.e., they do not contain the same number of up-pointing and down-pointing unit triangles), they have no lozenge tilings. This completes the proof. \epf


\proclaim{Proposition \tdb} $(${\rm a}$)$. Let $H_{x,y,z}(k,l)$ be the region obtained from the hexagon of side lengths $x$, $y+k+1$, $z$, $x+k+1$, $y$, $z+k+1$ $($clockwise from top$)$ by removing an up-pointing unit triangle from its northwestern side, $l$ units above the western corner, and an up-pointing triangle of side $k$ from its northeastern side, one unit above the eastern corner $($see Figure ${\fdc}$ for an illustration$)$.

\topinsert
\twoline{\mypic{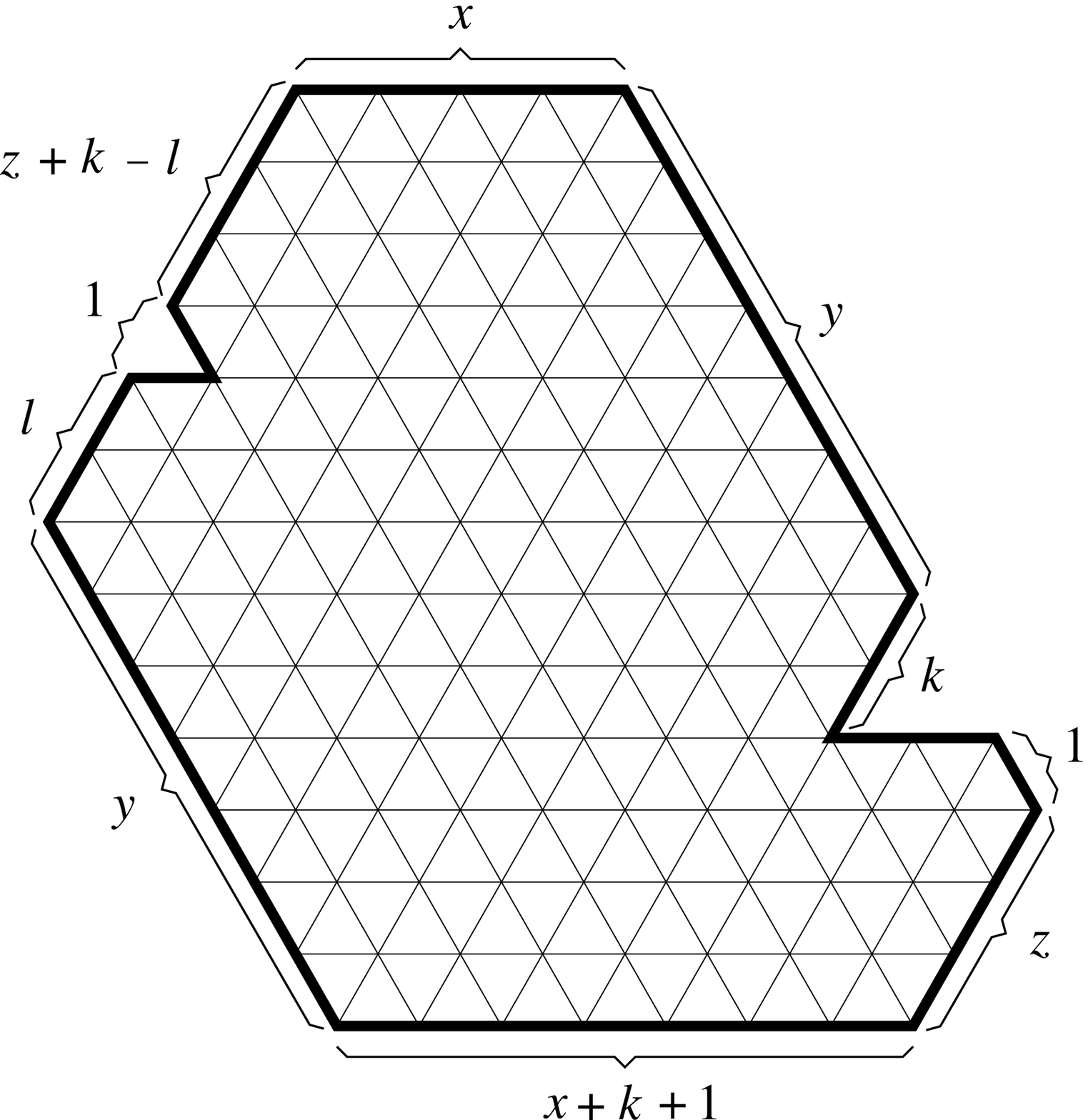}}{\mypic{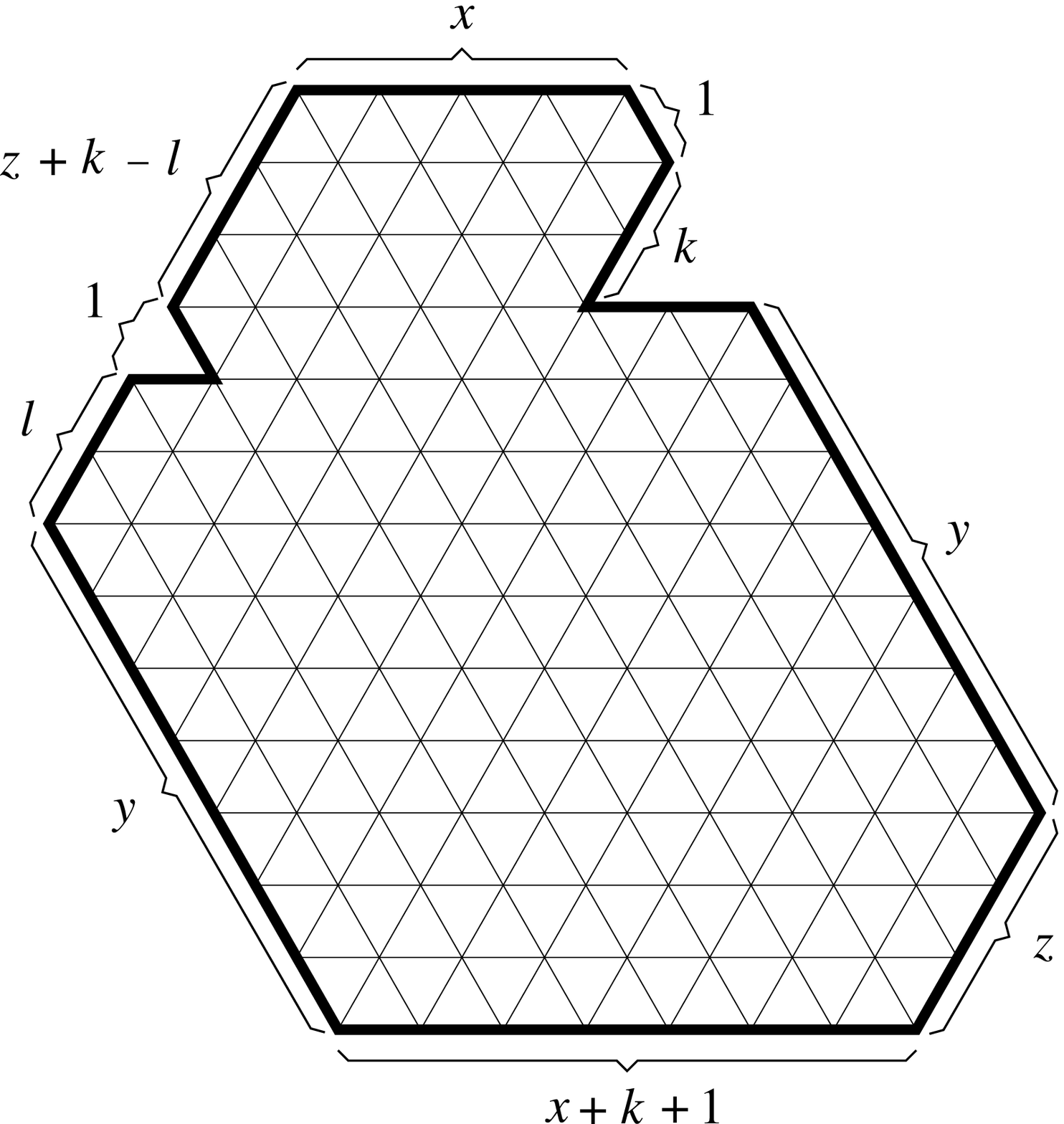}}
\medskip
\centerline{{\smc Figure~{\fdc}. {\rm The hexagons with two notches $H_{4,7,3}(2,2)$ (left) and $H'_{4,7,3}(2,2)$ (right).}}}
\endinsert

Let $m=\min(x,y)$ and $M=\max(x,y)$. Then we have
$$
\M(H_{x,y,z}(k,l))=\M(H_{x,y,k})\frac{p(z,l)}{p(0,0)},
\tag\edc
$$
where $\M(H_{x,y,k})$ is given by $(\ece)$, and the polynomial $p(z,l)$ is defined to be
$$
\spreadlines{2\jot}
\align
&
p(z,l):=(l+1)_y (z+k-l+1)_x 
\\
&\ \ \ \ 
\times
(z+k+2)(z+k+2)^2\cdots(z+k+m+1)^m (z+k+m+2)^m\cdots(z+k+M+1)^m 
\\
&\ \ \ \ 
\times
(z+k+M+2)^{m-1} (z+k+M+3)^{m-2}\cdots (z+k+M+m)
\\
&\ \ \ \ 
\times 
\sum_{i=1}^{k+1} \frac{(-1)^{i-1}}{(i-1)!(k-i+1)!} (l-k+i)_{k-i+1} (l+y+1)_{i-1} (z+1)_{i-1} (z+i+1)_{k-i+1}.
\tag\edd
\endalign
$$

$(${\rm b}$)$. Let $H'_{x,y,z}(k,l)$ be the region defined precisely as $H_{x,y,z}(k,l)$, with the one exception that the up-pointing triangle of side $k$ is one unit below the northeastern corner, rather than one unit above the eastern corner $($see Figure ${\fdc}$ for an illustration$)$.

Let $\nu=\min(y-1,k)$, and define $d(z)$ by
$$
\spreadmatrixlines{2\jot}
\align
\!\!\!\!\!\!\!\!\!\!\!\!
d(z):=\left\{\matrix (z+2)^1\cdots(z+\nu+1)^{\nu}\cdots(z+y+k-\nu)^{\nu}\cdots(z+y+k-1)^1, & \nu\geq1\\
1, & \nu=0\\
\dfrac{1}{(z+1)_k}, & \nu=-1\\
\endmatrix\right.
\tag\ede
\endalign
$$
$($in the first branch the bases are incremented by 1 from each factor to the next; the exponents are incremented by one until they reach $\nu$, stay equal to $\nu$ across the middle portion, and then they decrease by one unit from each factor to the next$)$.

Then we have
$$
\M(H'_{x,y,z}(k,l))={x+k\choose k}\frac{q(z,l)}{q(0,0)},
\tag\edf
$$
where the polynomial $q(z,l)$ is defined to be
$$
\spreadlines{2\jot}
\align
&
q(z,l):=d(z)\,(l+1)_y (z+k-l+1)_x 
\\
&\ \ \ \ 
\times
(z+k+2)(z+k+2)^2\cdots(z+k+m+1)^m (z+k+m+2)^m\cdots(z+k+M+1)^m 
\\
&\ \ \ \ 
\times
(z+k+M+2)^{m-1} (z+k+M+3)^{m-2}\cdots (z+k+M+m)
\\
&\, 
\times 
\sum_{i=1}^{k+1} \frac{(-1)^{i-1}}{(i-1)!(k-i+1)!} (l-k+i)_{k-i+1} (l+y+1)_{i-1} (l-k-z)_{i-1} (l-k-z+i)_{k-i+1}
\tag\edg
\endalign
$$
$($as in part \text{\rm (a)}, $m=\min(x,y)$ and $M=\max(x,y)$$)$.

\endproclaim

Note that the formulas giving $p(z,l)$ and $q(z,l)$ are very closely related: except for the factor $d(z)$ in the latter, the linear parts are precisely the same, and the sum factor in the latter is obtained from the sum factor in the former by replacing $z$ by $l-k-1-z$.

Furthermore, the constant multiple ${x+k\choose k}$ in the formula for $\M(H'_{x,y,z}(k,l))$ arises in fact as $\M(H_{x,1,k})$ (the two are equal by (\ece)), and is thus analogous to the constant multiple in the formula for $\M(H_{x,y,z}(k,l))$.

\medskip

Our proof of the above result is based on Kuo's original graphical condensation recurrence (see \cite{\KuoOne}). For ease of reference, we state below the particular instance of Kuo's general results that we need for our proofs (which is Theorem~2.4 in \cite{\KuoOne}).

\proclaim{Theorem {\tdc} (Kuo)} Let $G=(V_1,V_2,E)$ be a plane bipartite graph in which $|V_1|=|V_2|+1$. Let vertices $a$, $b$, $c$ and $d$ appear cyclically on a face of $G$. If $a,b,c\in V_1$ and $d\in V_2$, then
$$
\M(G-b)\M(G-\{a,c,d\})=\M(G-a)\M(G-\{b,c,d\})+\M(G-c)\M(G-\{a,b,d\}).
\tag\edh
$$

\endproclaim

 {\it Proof of Proposition {\tdb}.}\ We prove part (a) by induction, using Kuo condensation. Augment the region $H_{x,y,z}(k,l)$ by placing on its top a trapezoidal band consisting of $2x-1$ unit triangles as illustrated in Figure {\fdd} (which shows the region obtained this way from the region on the left in Figure {\fdc}); denote the resulting region by $\tilde{H}_{x,y,z}(k,l)$. We apply Kuo condensation to the dual graph of $\tilde{H}_{x,y,z}(k,l)$, with the vertices $a,b,c,d$ corresponding to the unit triangles indicated in Figure {\fdd}.

\topinsert
\twoline{\mypic{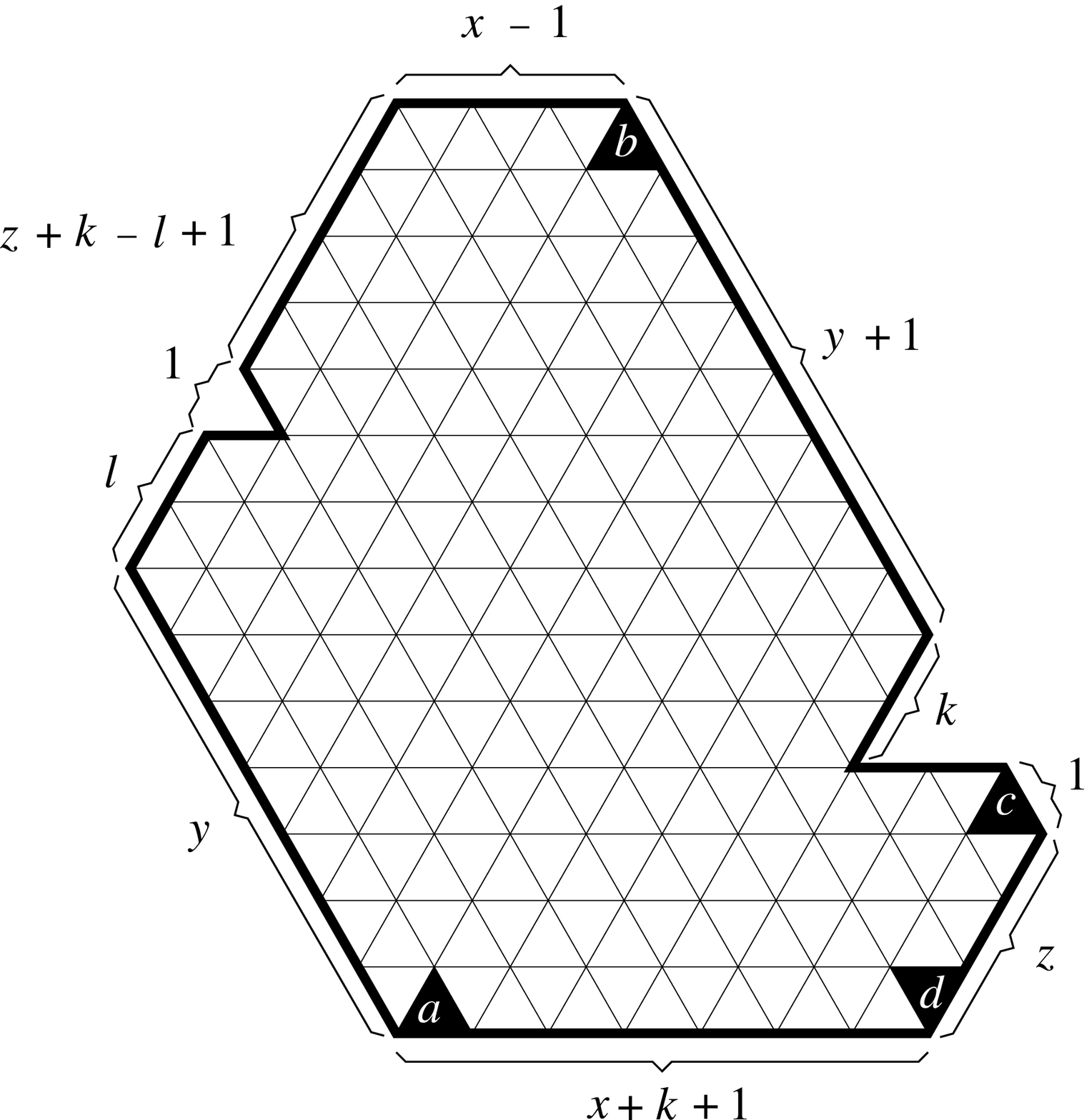}}{\mypic{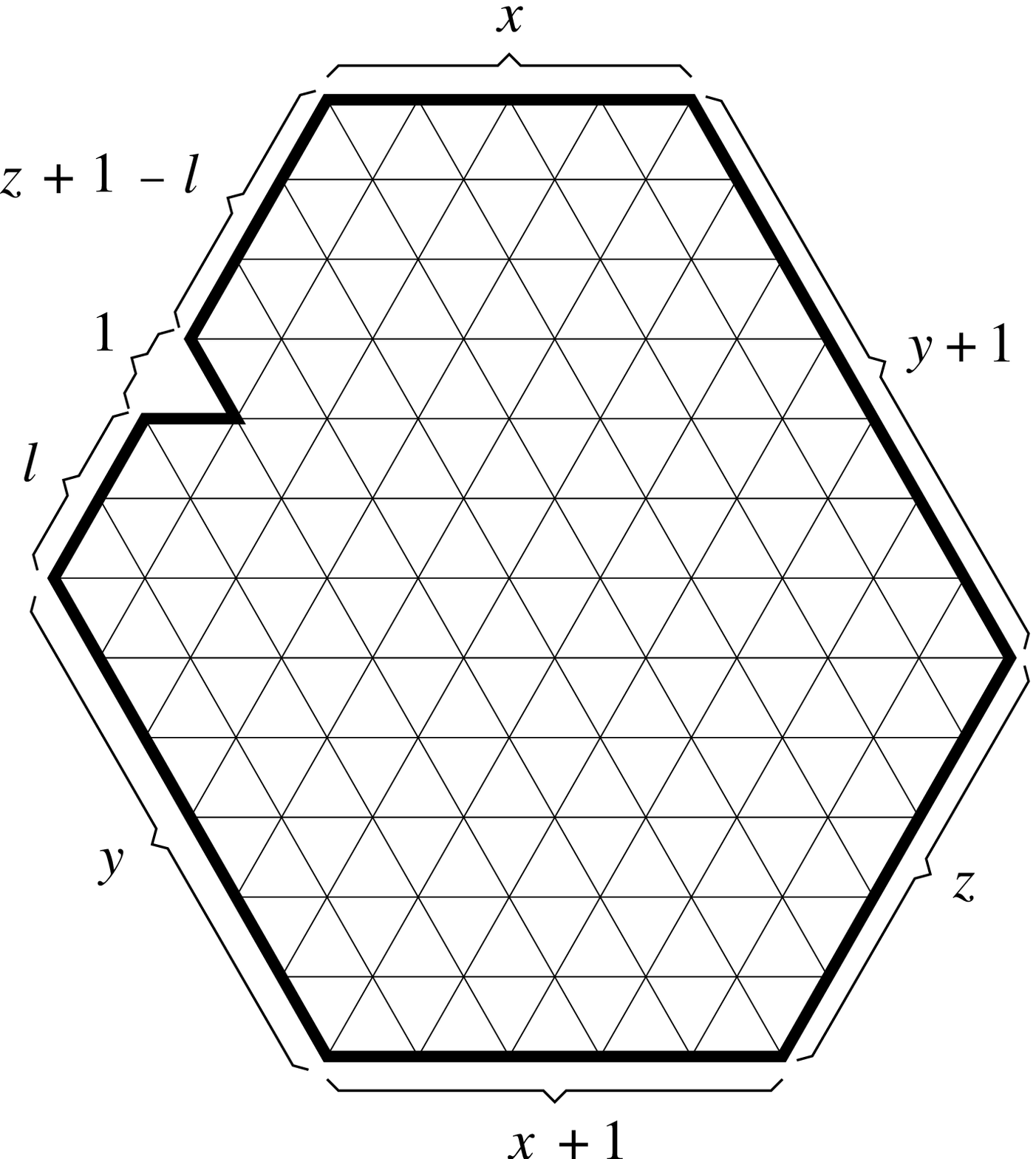}}
\medskip
\twoline{{\smc Figure~{\fdd}. {\rm The augmented region $\tilde{H}_{4,7,3}(2,2)$.}}}{{\smc Figure~{\fddp}. {\rm $F_{4,6,5}(2)$.}}}
\endinsert

\topinsert
\twoline{\mypic{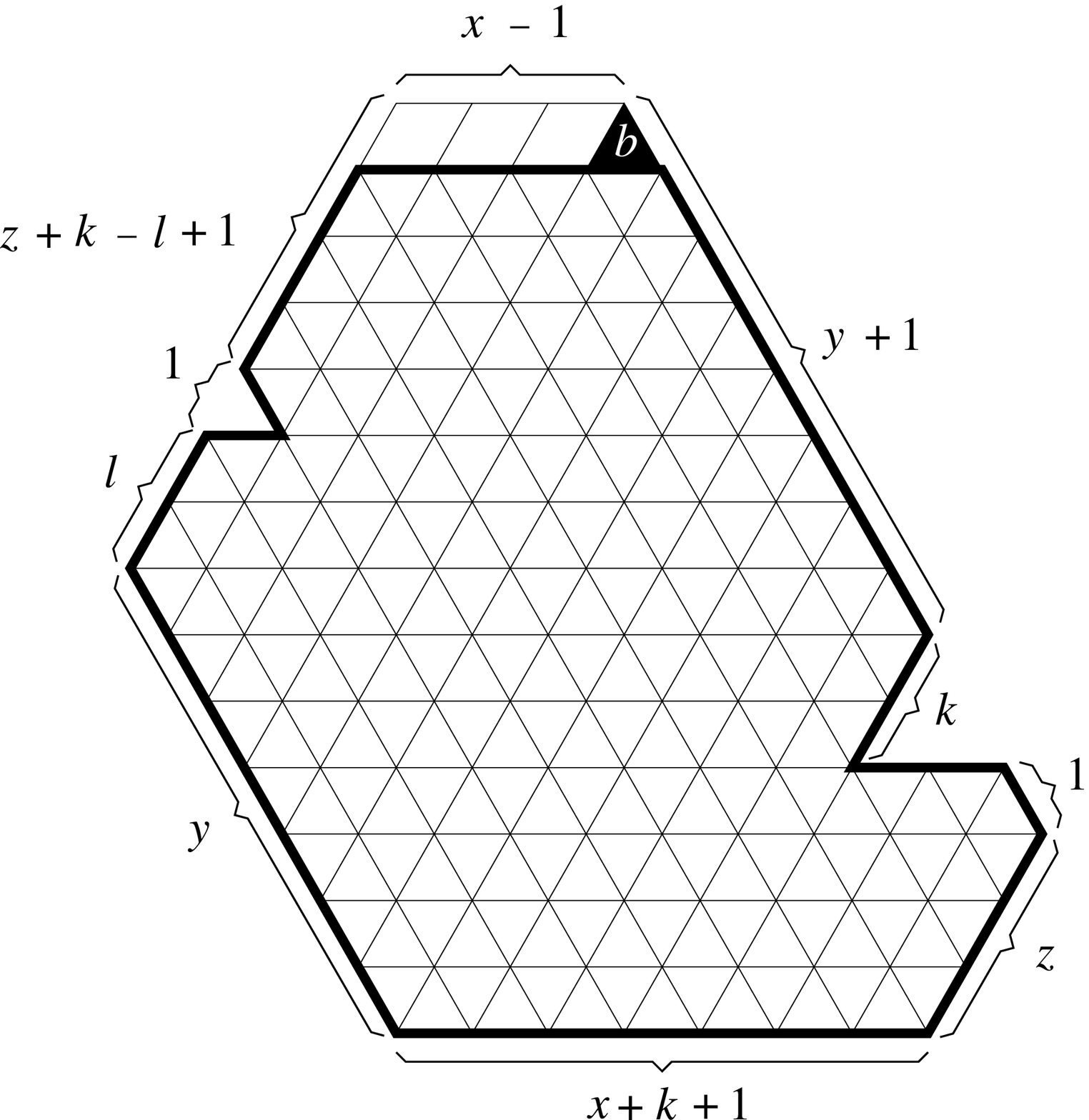}}{\mypic{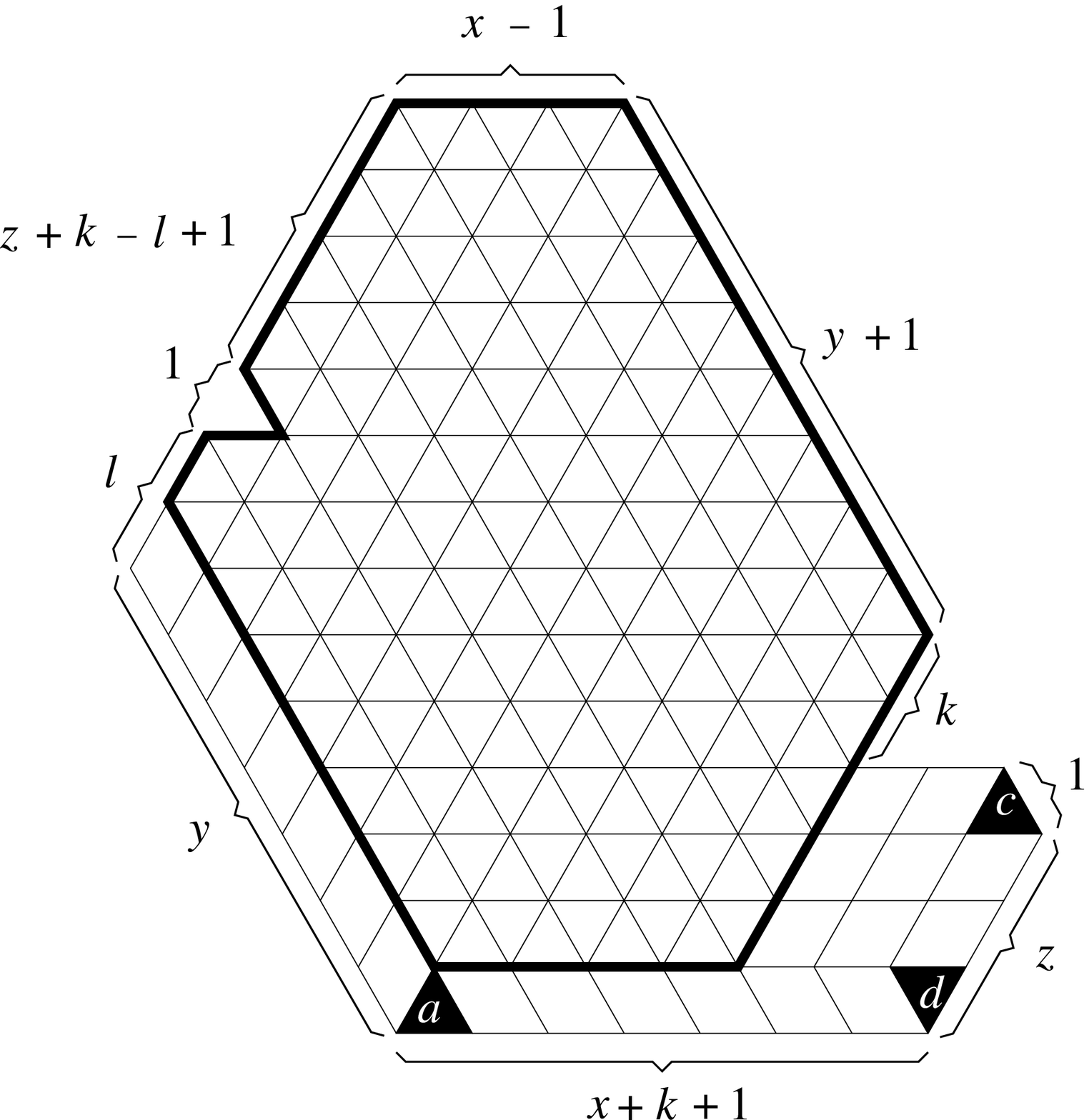}}
\medskip
\twoline{\mypic{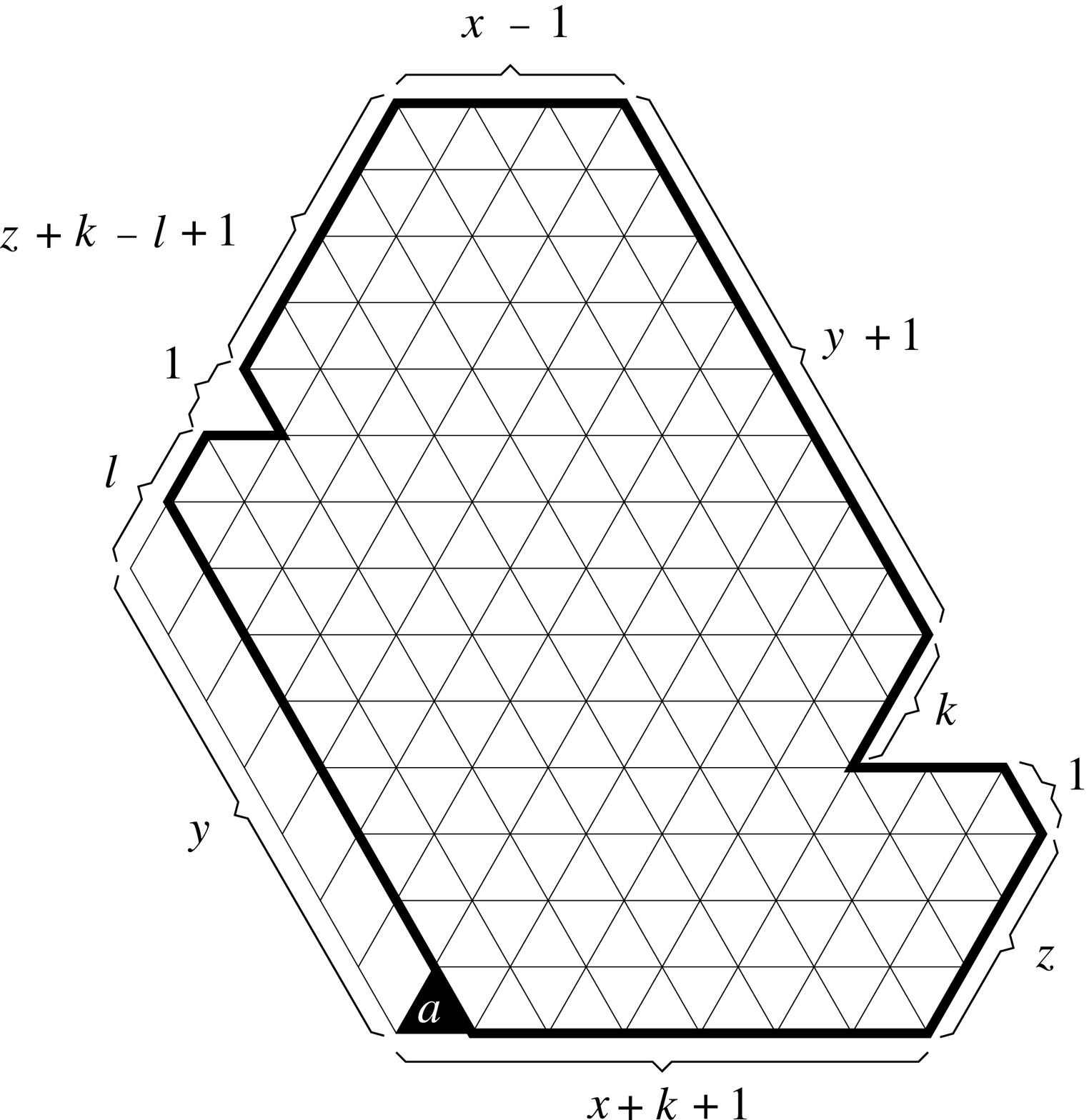}}{\mypic{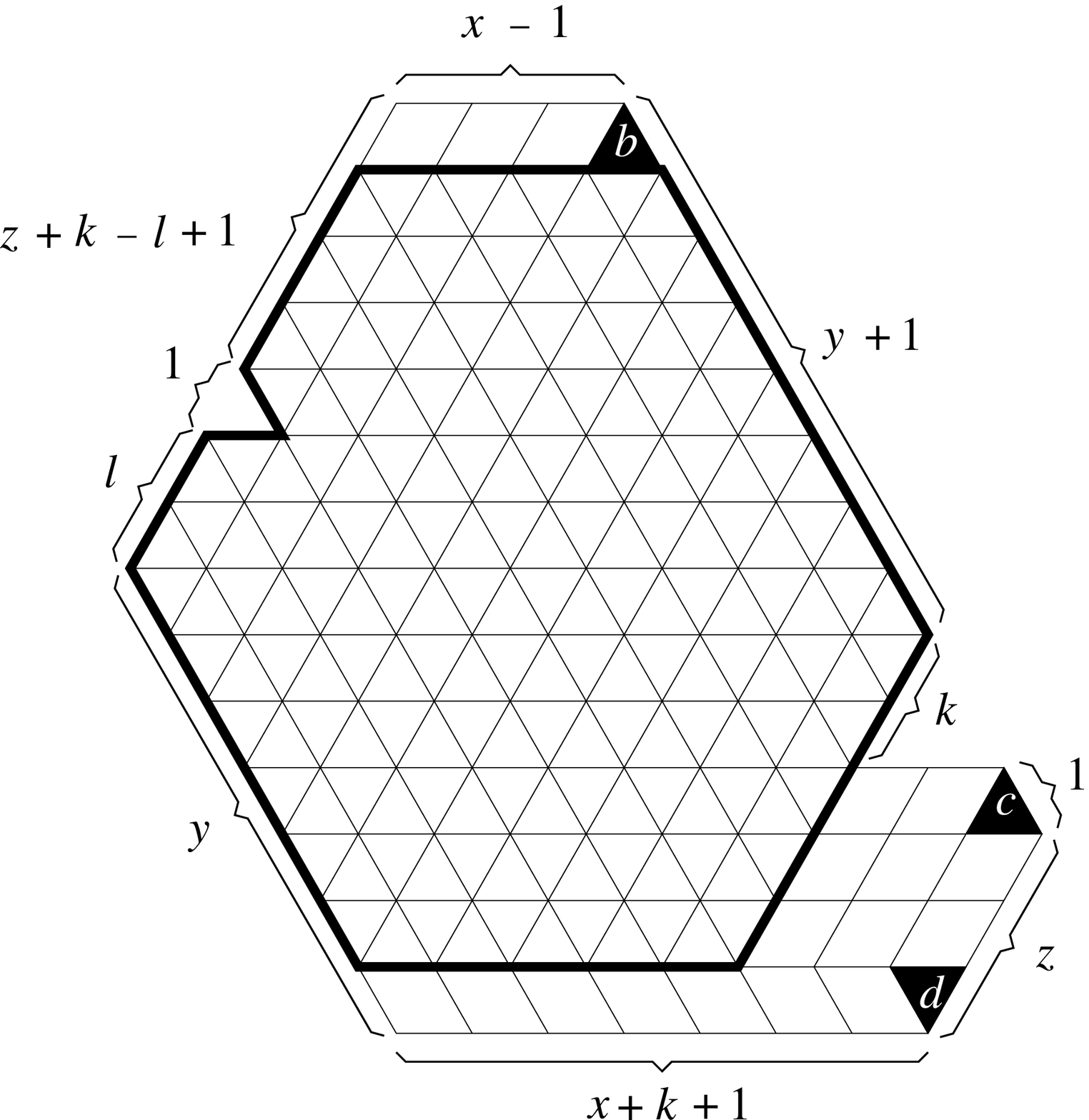}}
\medskip
\twoline{\mypic{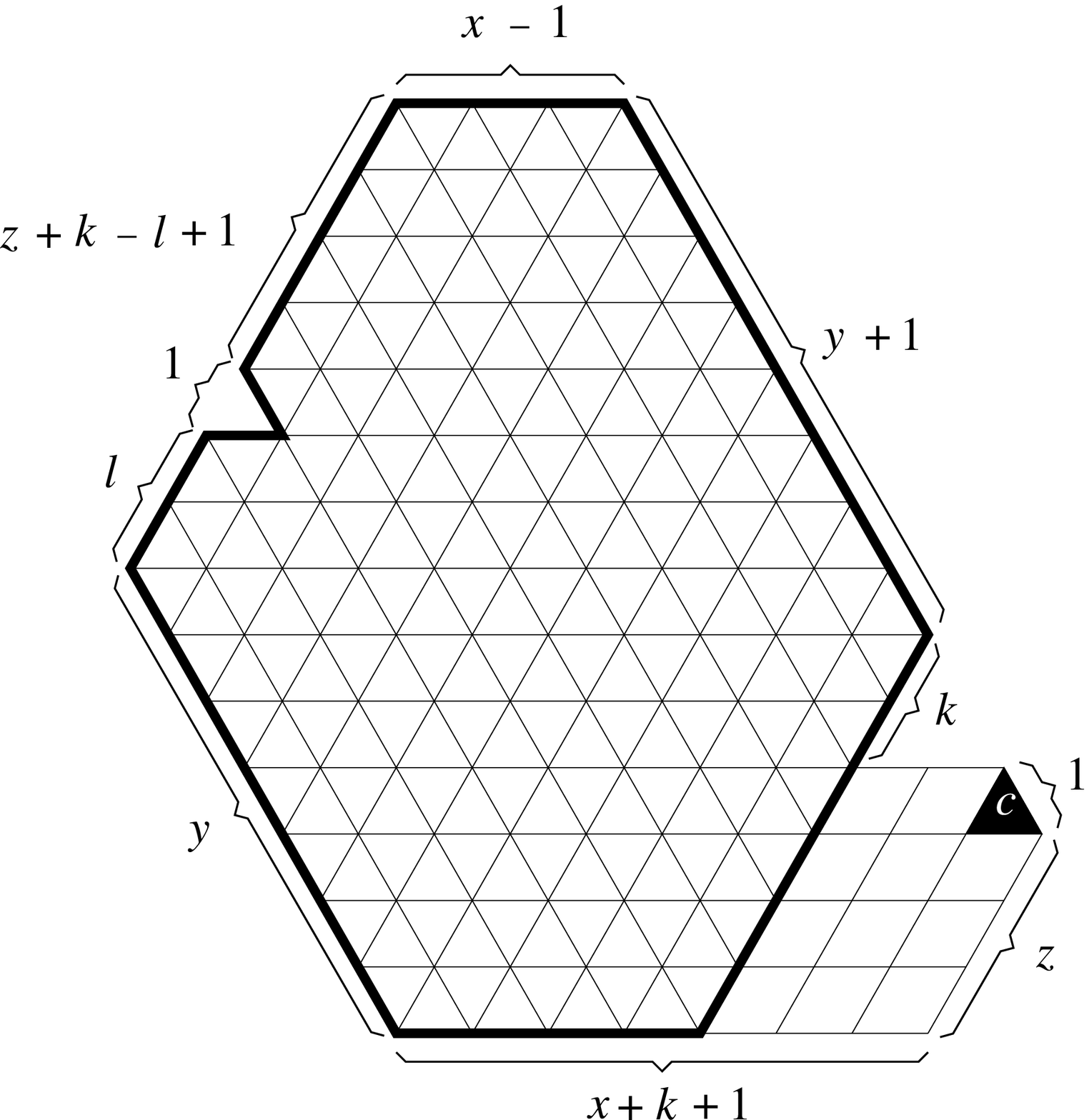}}{\mypic{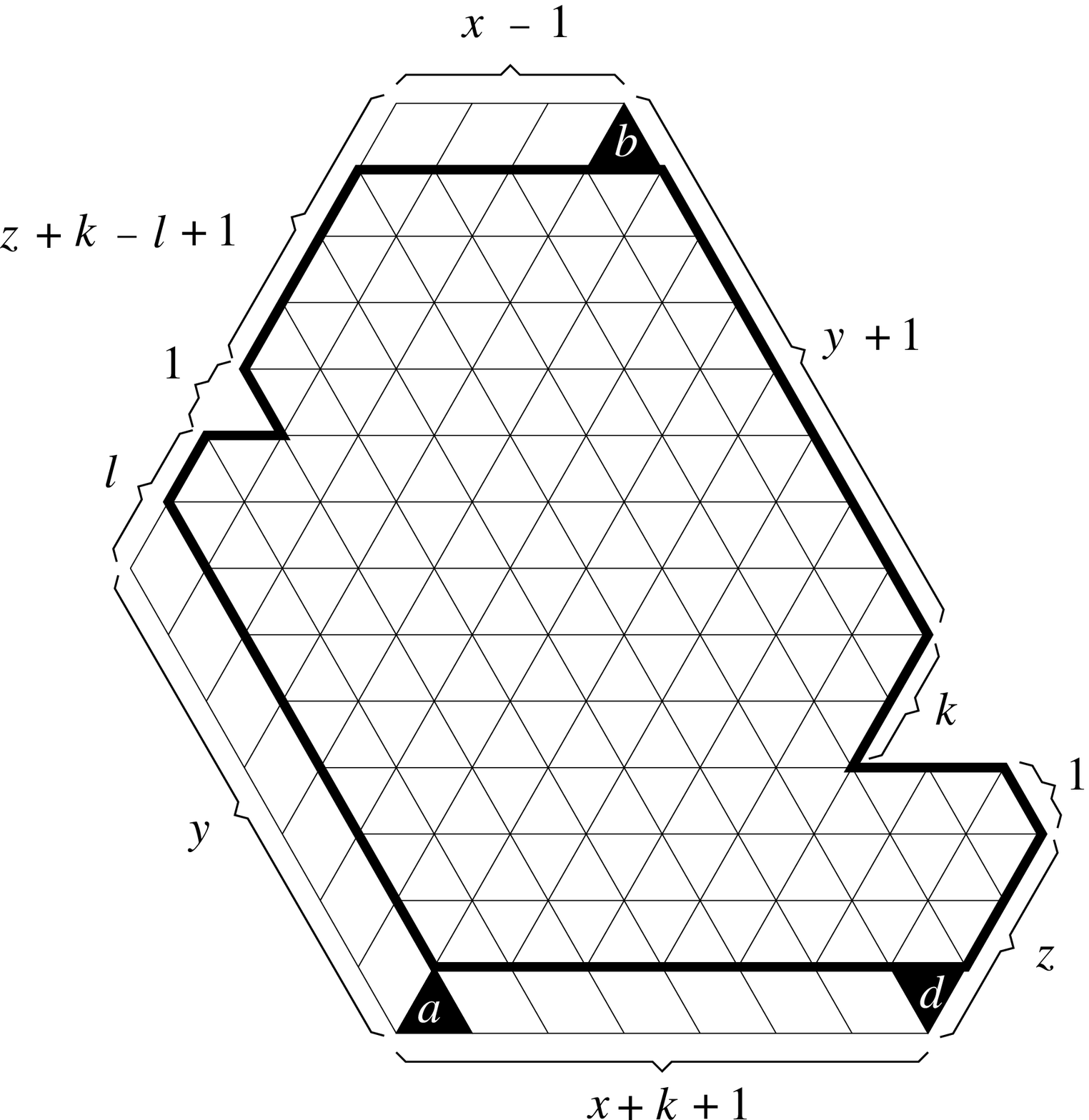}}
\medskip
\centerline{{\smc Figure~{\fde}. {\rm Obtaining the recurrence for the regions $H_{x,y,z}(k,l)$.}}}
\endinsert

In the graphs resulting this way from equation (\edh) there are many edges that are forced to be part of every perfect matching. The situation is illustrated  --- for the corresponding dual lattice regions --- in Figure {\fde}. 

Equation (\edh) states that the product of the number of lozenge tilings of the two regions on top in Figure {\fde} is equal to the product of the number of lozenge tilings of the two regions in the middle, plus the product of the number of lozenge tilings of the two regions on the bottom. 

After removing the forced lozenges, it is clear that the region resulting from the region on the top left in  Figure {\fde} is $H_{x,y,z}(k,l)$. 

Denote by $F_{x,y,z}(l)$ the region obtained from a hexagon of sides $x,y+1,z,x+1,y,z+1$ by removing the up-pointing unit triangle from its boundary that is $l$ units above the western corner (see Figure {\fddp}). Then what is left from the region on the top right in  Figure {\fde} after removing the forced lozenges is precisely the region $F_{x-1,y,z+k}(l-1)$.

Similarly, one sees that the regions resulting from the two regions in the middle of Figure~{\fde} after removing the forced lozenges are $H_{x-1,y+1,z}(k,l-1)$and $F_{x,y-1,z+k}(l)$. The two regions on the bottom in Figure {\fde} lead similarly to $F_{x-1,y,z+k+1}(l)$ and $H_{x,y,z-1}(k,l-1)$, respectively.

Therefore, by equation (\edh) we obtain
$$
\spreadlines{2\jot}
\align
&
\M(H_{x,y,z}(k,l)\M(F_{x-1,y,z+k}(l-1))
=
\M(H_{x-1,y+1,z}(k,l-1))\M(F_{x,y-1,z+k}(l))
\\
&\ \ \ \ \ \ \ \ \ \ \ \ \ \ \ \ \ \ \ \ \ \ \ \ \ \ \ \ \ \ \ \ 
+
\M(H_{x,y,z-1}(k,l-1))\M(F_{x-1,y,z+k+1}(l)).
\tag\edi
\endalign
$$
All the regions in this recurrence are well-defined provided $x\geq1$, $z\geq1$, and $l\geq1$.

\topinsert
\threeline{\mypic{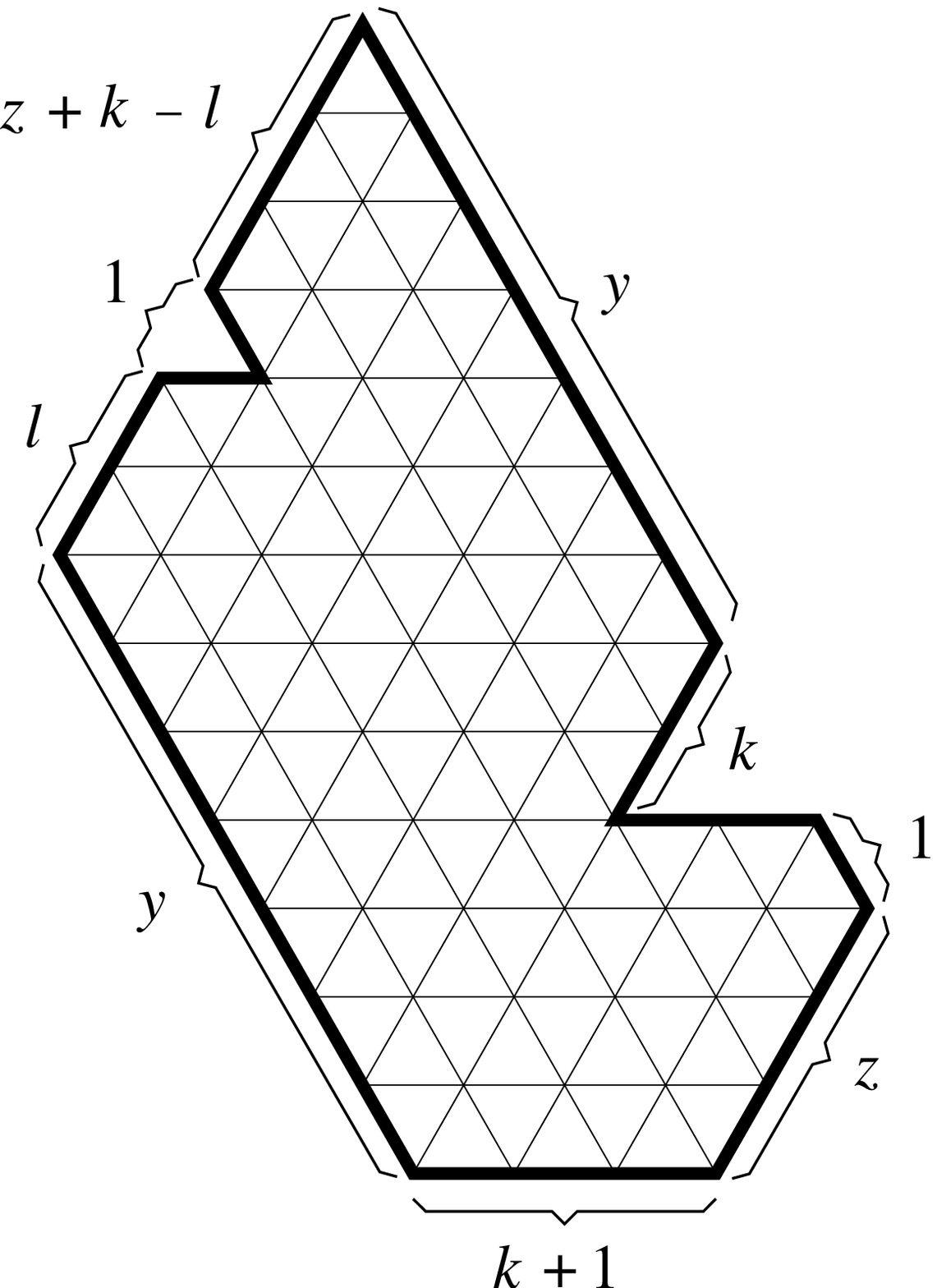}}{\mypic{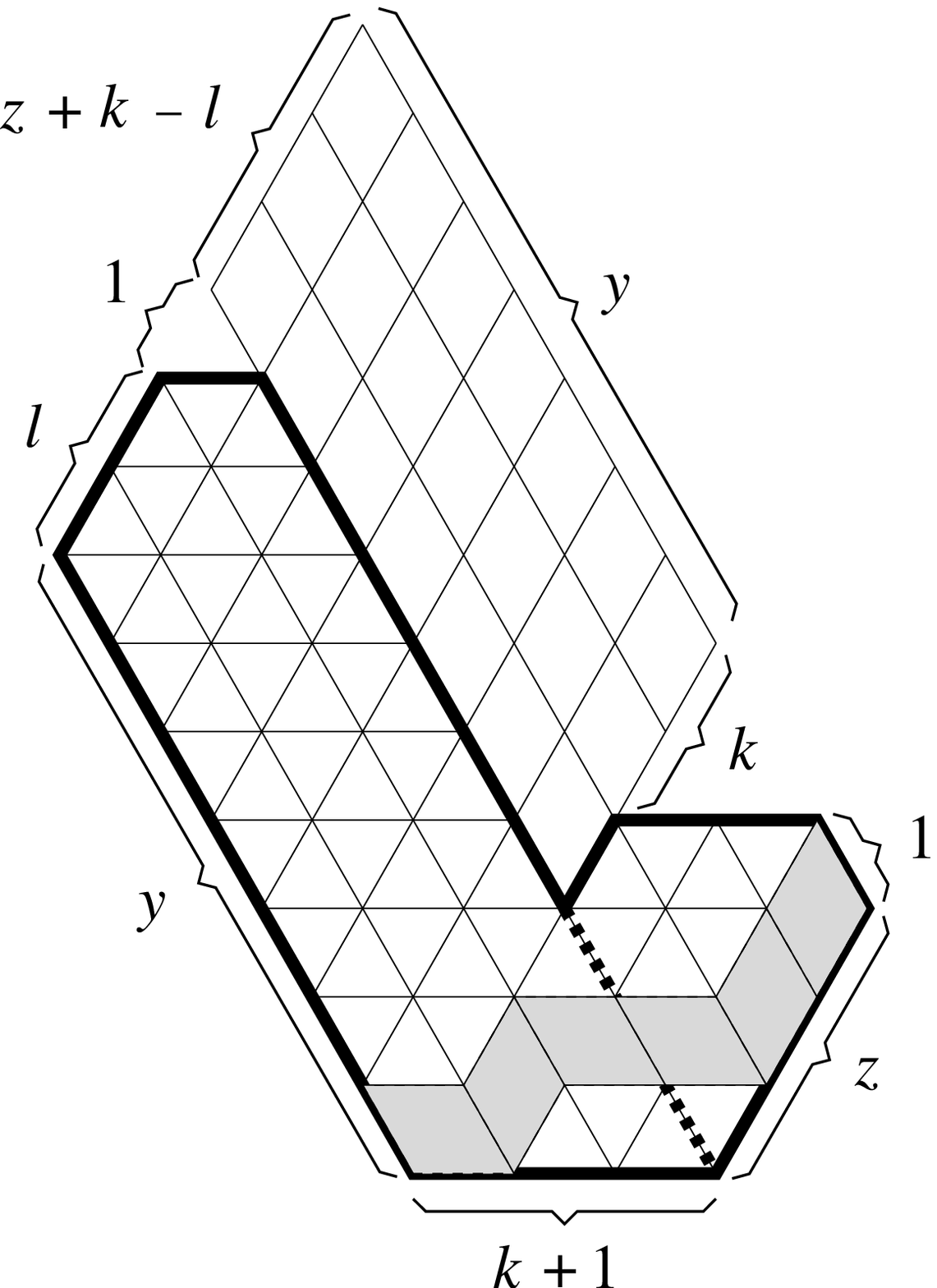}}{\mypic{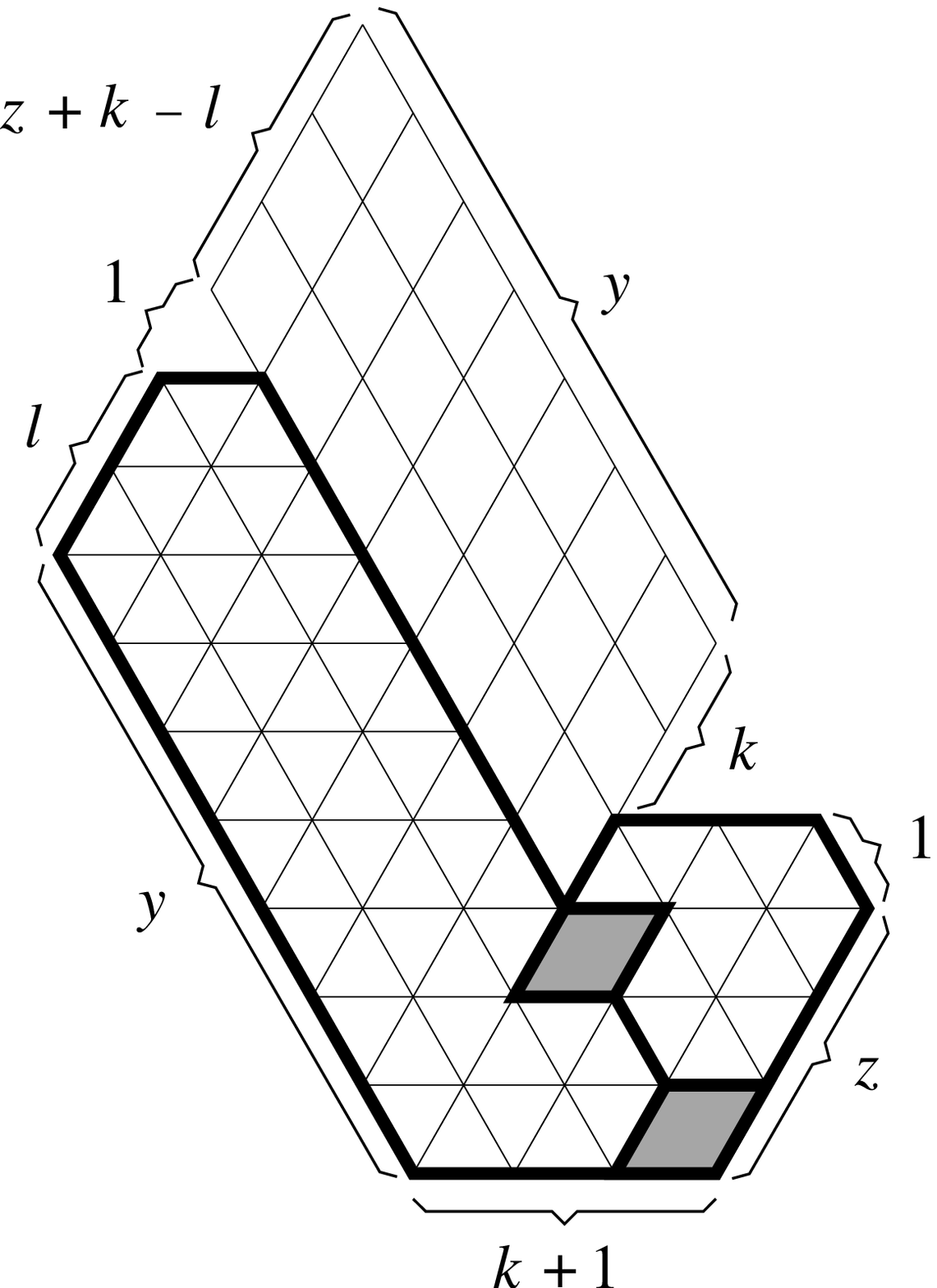}}
\medskip
\threeline{\text{\rm \ \ \ \ \ \ \ (a)\ \ \ \ \ \ \ \ \ \ \ \ \ \ \ \ \ \ \ \ \ }}{\text{\rm (b)}}{\text{\rm \ \ \ \ \ \ \ \ \ \ \ \ \ \ \ \ \ \ \ \ \ (c)}}
\medskip
\centerline{{\smc Figure~{\fddpp}. {\rm The case $x=0$.}}}
\endinsert

\topinsert
\twoline{\mypic{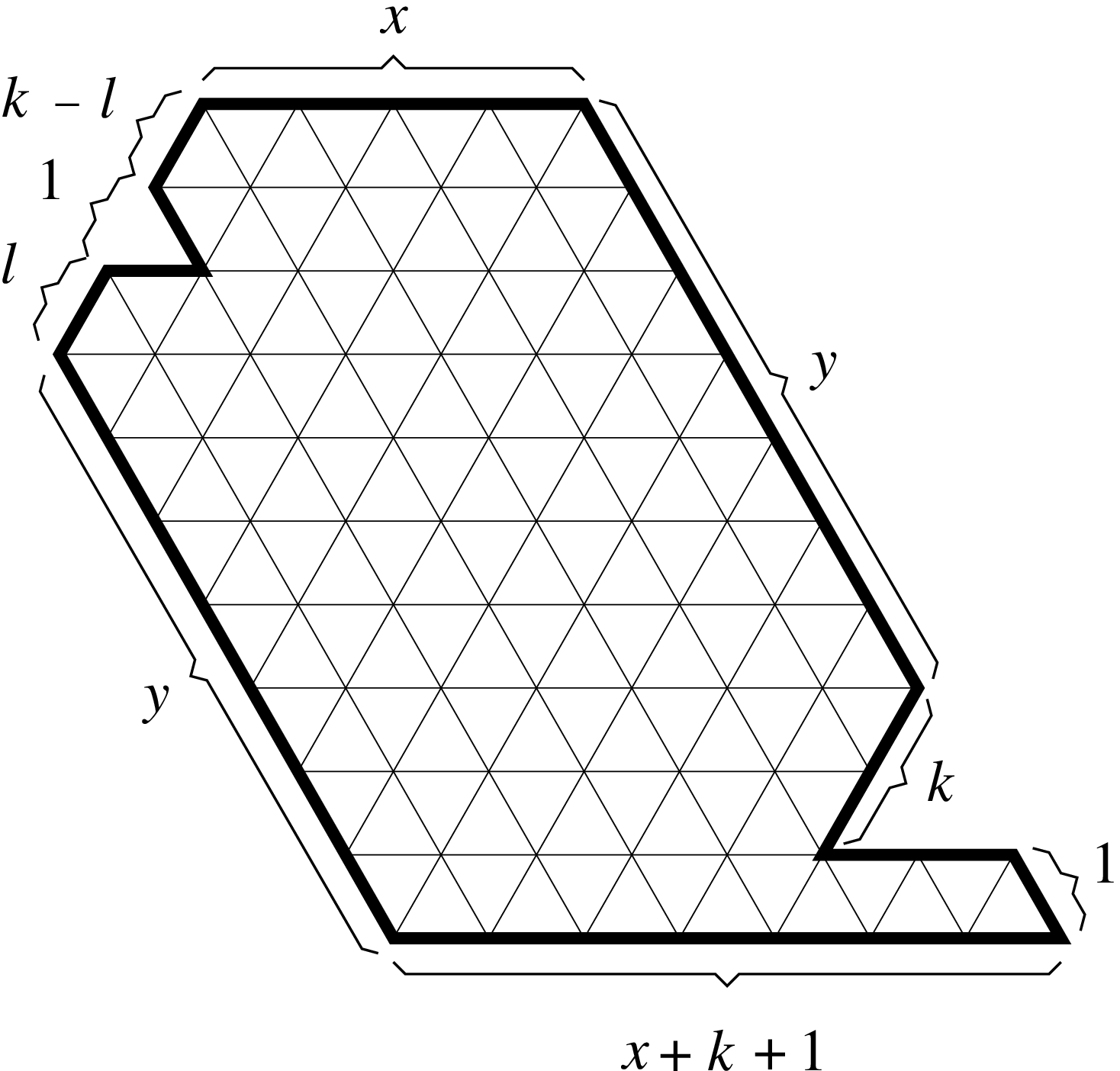}}{\mypic{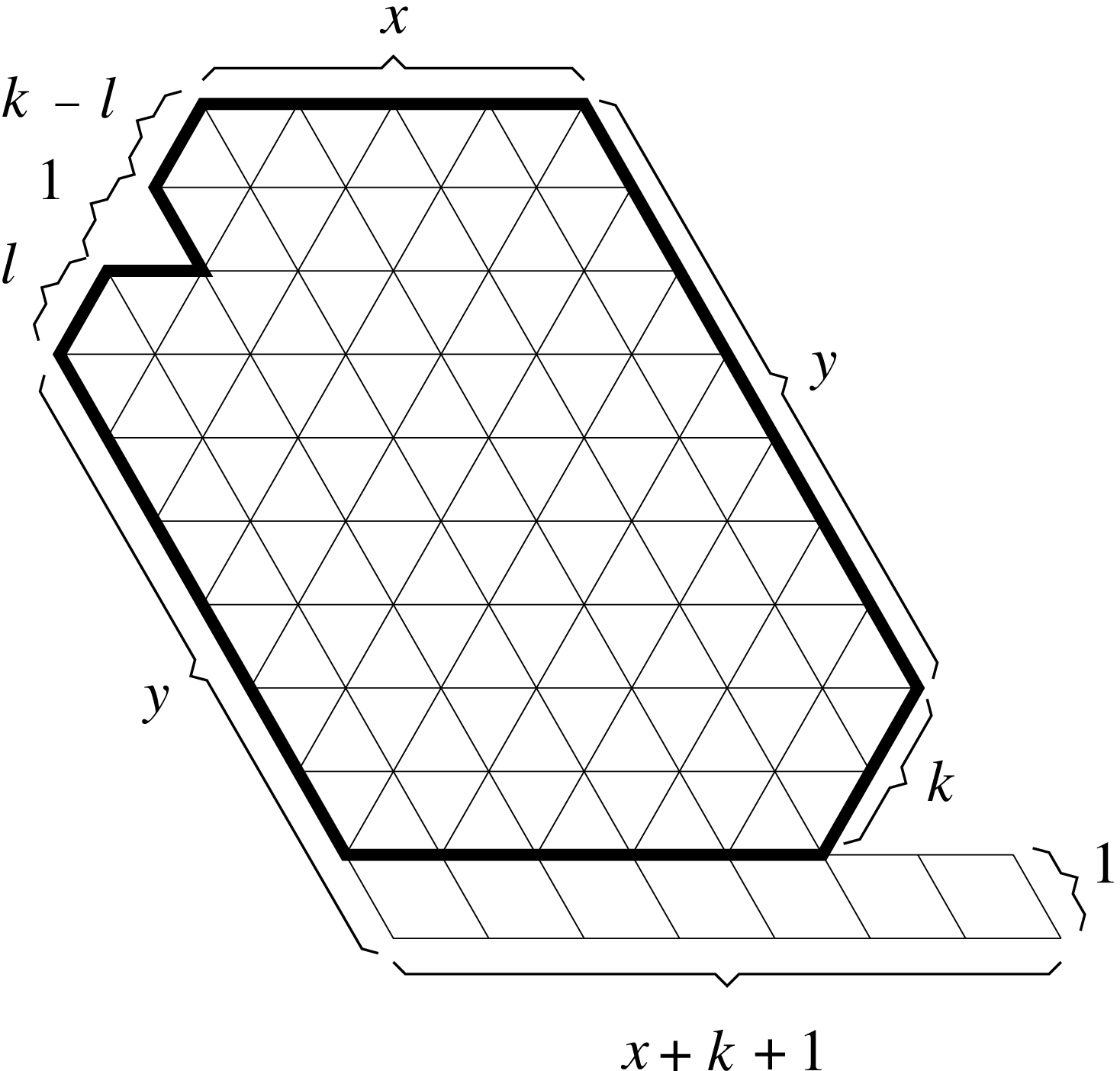}}
\medskip
\twoline{\text{\rm \ \ \ \ \ \ \ (a)\ \ \ \ \ \ \ \ \ \ \ \ \ \ \ \ \ \ \ \ \ }}{\text{\rm (b)}}
\medskip
\centerline{{\smc Figure~{\fddppp}. {\rm The case $z=0$.}}}
\endinsert

For $x=0$ we can verify formula (\edc) directly. Indeed, in this case the region $H_{0,y,z}(k,l)$ looks as illustrated in Figure {\fddpp}(a). Note that after removing the forced lozenges, the path of lozenges connecting the portion of length 1 of the boundary just above the eastern corner to the southwestern side (which necessarily ends at the bottommost unit segment of the latter) must cross the thick dotted line in Figure {\fddp}(b). This dotted line has length $k+1$ (or possibly less, in case it meets the boundary at an interior point of the southern side). The number of tilings of $H_{0,y,z}(k,l)$ for which the segment $s$ at which this path of lozenges meets the dotted line is fixed is readily seen to be the same as the number of tilings of $H_{0,y,z}(k,l)$ which contain all the lozenges that straddle the remaining unit segments of the dotted line, but not the lozenge that straddles the unit segment $s$ (see Figure {\fddpp}(c)). However, the latter is just the product of the number of tilings of the two $T$-type regions indicated in Figure {\fddpp}(c). Therefore, by Proposition {\tcb}, $\M(H_{0,y,z}(k,l))$ is equal to a sum of (at most) $k+1$ simple products. It is straightforward to check that the resulting expression agrees with the right hand side of (\edc).

For $z=0$, we can verify formula (\edc) as follows. The region $H_{x,y,0}(k,l)$ looks as illustrated in Figure {\fddppp}(a). After removing the forced lozenges,we obtain from it a $T$-type region, whose number of tilings is given by Proposition {\tcb}. It is not hard to check that the resulting formula agrees with the $z=0$ specialization of the expression on the right hand side of (\edc).
Indeed, by construction $z+k-l\geq0$ (see the picture on the left in Figure {\fdc}), and since we are in the case $z=0$, we have $l\leq k$. 
If $l<k$, due to the factor $(l-k+i)_{k-i+1}$ in the summand of the sum in (\edd), this summand is non-0 only for $i=k+1$. Thus verification of (\edc) amounts to checking that two explicit products of linear factors have the same value, which is readily checked. 

For $l=k$, the sum in (\edd) can be written in terms of hypergeometric series\footnote{The hypergeometric function of parameters
$a_1,\dotsc,a_p$ and $b_1,\dotsc,b_q$ is defined by
$${}_p F_q\!\left[\matrix a_1,\dotsc,a_p\\ b_1,\dotsc,b_q\endmatrix;
z\right]=\sum _{k=0} ^{\infty}\frac {(a_1)_k\cdots(a_p)_k}
{k!\,(b_1)_k\cdots(b_q)_k} z^k. $$} as 
$$
\spreadlines{3\jot}
\align
&
\sum_{i=1}^{k+1} \frac{(-1)^{i-1}}{(i-1)!(k-i+1)!} (i)_{k-i+1} (k+y+1)_{i-1} (z+1)_{i-1} (z+i+1)_{k-i+1}
=
\\
&\ \ \ \ \ \ \ \ \ \ \ \ \ \ \ \ \ \ \ \ \ \ \ \ \ \ \ \ \ \ \ \ \ \ \ \ \ \ \ \ \ \ \ \ \ \ \ \ \ \ \ \ \ \ \ \  
{}_3F_2\!\left[\matrix{-k,\,z+1,\,k+y+1}\\{1,\,z+2}\endmatrix; 1\right](z+2)_k.
\endalign
$$
Since $z=0$ in the case under consideration, the right hand side above becomes
$$
(k+1)!\,\,
{}_2F_1\!\left[\matrix{-k,\,k+y+1}\\{2}\endmatrix; 1\right],
$$
which by Lemma {\tdd} evaluates to an explicit product of linear factors.
Thus the verification of case $z=0$, $l=k$ also amounts to checking that two explicit products of linear factors agree, which is easily done. 
%
%

We may assume therefore that $x\geq1$ and $z\geq1$. We prove formula (\edc) by induction on $l$, using recurrence (\edi) at the induction step.

The base case is $l=0$. It is clear from Figure {\fdc} that for $l=0$ there is a band of forced lozenges along the southwestern side of $H_{x,y,z}(k,l)$, and that after removing this band one is left with a hexagon with a single notch of side $k$ on its northeastern side. However, such a region is readily seen to be a $T$-region of the type addressed by Proposition {\tcb} (see  Figure {\fcg}), and thus the number of its lozenge tilings is given by the product formula (\ech). It is routine to check that the resulting formula agrees with the $l=0$ specialization of formula (\edc) (note in particular that, due to the presence of the factor $(l-k+i)_{k-i+1}$ in the summand in (\edd), all but the last term in the sum in (\edd) are zero).

For the induction step, assume that formula (\edc) holds for all instances when the value of the $l$-parameter is $l-1$, and consider the region $H_{x,y,z}(k,l)$. Since we are in the case $x\geq1$ and $z\geq1$, and we are assuming $l\geq1$, all six regions in equation (\edi) are well defined. Moreover, the $F$-regions are clearly special cases of the $T$-regions addressed by  Proposition {\tcb}, and therefore have the number of their lozenge tilings expressed by the simple product formula (\ech). 

Therefore, by (\edi) and by the induction hypothesis, we obtain that $\M(H_{x,y,z}(k,l))$ is a sum of two concrete product expressions, each involving, besides a single factor having the type of the sum in (\edd), only linear factors. It is routine to check that the sum of these two products agrees with the product on the right hand side of (\edc). This concludes the proof of part (a).

The proof of part (b) is completely analogous.
\epf

\proclaim{Lemma \tdd} For any non-negative integers $k$, $y$ and $z$ we have
$$
{}_2F_1\!\left[\matrix{-k,\,k+y}\\{z}\endmatrix; 1\right]
=
\frac{(z-y-k)_k}{(z)_k}.
\tag\edj
$$

\endproclaim

\pf  By Gauss' formula (see e.g. \cite{\Sl, (1.7.6),\,Appendix (III.3)}), for any $a,b,c\in\C$ with $\Rep\,(c-a-b)>0$ and $c\neq0,-1,-2,\dotsc,$ 
one has
$$
{}_2 F_1\left[\matrix a,b\\ c\endmatrix;1\right]=
\frac{\Gamma(c)\Gamma(c-a-b)}{\Gamma(c-a)\Gamma(c-b)}.
$$
This formula applies to the hypergeometric series in the statement of the lemma provided $z\geq y+1$, and proves the statement in this case. Note that, due to the presence of the numerator parameter $-k$, the sum on the left hand side of (\edj) is finite. Furthermore, after multiplication by $(z)_k$, both sides of (\edj) become polynomials in $z$. Since we have seen that the identity holds for infinitely many values of $z$, it follows that it holds for all~$z$.~\epf


\mysec{References}
{\openup 1\jot \frenchspacing\raggedbottom
\roster

\myref{\anglep}
  M. Ciucu and I. Fischer, A triangular gap of side 2 in a sea of dimers in a $60^\circ$ angle, {\it J. Phys. A: Math. Theor.} {\bf 45} (2012), 494011.


\myref{\ff}
 M Ciucu and C. Krattenthaler, A dual of MacMahon's theorem on plane partitions, {\it Proc. Nat. Acad. Sci. USA} {\bf 110} (2013), 4518--4523.

\myref{\hexnotch}
  M. Ciucu and I. Fischer, Proof of two conjectures of Ciucu and Krattenthaler on the enumeration of lozenge tilings of hexagons with cut off corners, arxiv preprint arXiv:1309.4640, 2013.

\myref{\CLP}
  H. Cohn, M. Larsen, and J. Propp, The shape of a typical boxed plane partition,{\it New York J. of Math.} {\bf 4} (1998), 137--165.

\myref{\Eisen}
  T. Eisenk\"olbl, Rhombus Tilings of a Hexagon with Three Fixed Border Tiles, {\it J. Comb. Theory Ser. A}, {\bf 88} (1999), 368--378.

\myref{\Fulmek}
  M. Fulmek, Graphical condensation, overlapping Pfaffians and superpositions of matchings, {\it Electron. J. Combin.} {\bf 17} (2010), \#R83.

\myref{\GT}
  I. M. Gelfand and M. L. Tsetlin, Finite-dimensional representations of the group of unimodular matrices (in Russian), {\it Doklady Akad. Nauk. SSSR (N. S.)} {\bf 71} (1950), 825--828.

\myref{\KuoOne}
  E. H. Kuo, Applications of graphical condensation for enumerating matchings and tilings, {\it Theoret. Comput. Sci.} {\bf 319} (2004), 29--57.


\myref{\KuoTwo}
  E. H. Kuo, Graphical condensation generalizations involving Pfaffians and determinants, arxiv preprint math:CO/06055154, 2006.

\myref{\MacM}
  P. A. MacMahon, Memoir on the theory of the partition of numbers---Part V. 
Partitions in two-dimensional space, {\it Phil. Trans. R. S.}, 1911, A.

\myref{\Sl}
  L. J. Slater, ``Generalized hypergeometric series,'' Cambridge University Press, Cambridge, 1966.

\endroster\par}

\enddocument